\documentclass[10pt]{article}

\usepackage{pifont}
\usepackage{amsmath}
\usepackage{cases}
\usepackage{amsfonts}
\usepackage{latexsym}
\usepackage{amssymb}
\usepackage{stmaryrd}
\usepackage{overpic}
\usepackage{color}
\usepackage{abstract}
\usepackage{indentfirst}
\usepackage{booktabs}
\usepackage{url}
\usepackage{verbatim}
\usepackage{diagbox}
\usepackage{multirow} 
\usepackage{graphicx}  
\usepackage{float}  
\usepackage{subfig}

\DeclareFontEncoding{FMS}{}{}
\DeclareFontSubstitution{FMS}{futm}{m}{n}
\DeclareFontEncoding{FMX}{}{}
\DeclareFontSubstitution{FMX}{futm}{m}{n}
\DeclareSymbolFont{fouriersymbols}{FMS}{futm}{m}{n}
\DeclareSymbolFont{fourierlargesymbols}{FMX}{futm}{m}{n}
\DeclareMathDelimiter{\VERT}{\mathord}{fouriersymbols}{152}{fourierlargesymbols}{147}

\newtheorem{theorem}{Theorem}[section]

\newtheorem{example}[theorem]{Example}

\newtheorem{remark}[theorem]{Remark}

\begin{document}

\begin{center}
\large\bf Local Randomized Neural Networks Methods for Interface Problems
\end{center}

\vspace{2mm}

\begin{center}
 {\large\sc Yunlong Li}\footnote{School of Mathematics and Statistics, Xi'an Jiaotong University, Xi'an,
Shaanxi 710049, P.R. China. E-mail: {\tt 4122107033@stu.xjtu.edu.cn}},\quad {\rm and}\quad
{\large\sc Fei Wang}\footnote{School of Mathematics and Statistics, Xi'an Jiaotong University,
Xi'an, Shaanxi 710049, China. The work of this author was partially supported by
the National Natural Science Foundation of China (Grant No.\ 12171383). Email: {\tt feiwang.xjtu@xjtu.edu.cn}} 
\end{center}

\vspace{2mm}

\begin{quote} 
\noindent{}{\bf Abstract}: 
Accurate modeling of complex physical problems, such as fluid-structure interaction, requires multiphysics coupling across the interface, which often has intricate geometry and dynamic boundaries. Conventional numerical methods face challenges in handling interface conditions. Deep neural networks offer a mesh-free and flexible alternative, but they suffer from drawbacks such as time-consuming optimization and local optima. In this paper, we propose a mesh-free approach based on Randomized Neural Networks (RNNs), which avoid optimization solvers during training, making them more efficient than traditional deep neural networks. Our approach, called Local Randomized Neural Networks (LRNNs), uses different RNNs to approximate solutions in different subdomains. We discretize the interface problem into a linear system at randomly sampled points across the domain, boundary, and interface using a finite difference scheme, and then solve it by a least-square method. For time-dependent interface problems, we use a space-time approach based on LRNNs. We show the effectiveness and robustness of the LRNNs methods through numerical examples of elliptic and parabolic interface problems. We also demonstrate that our approach can handle high-dimension interface problems. Compared to conventional numerical methods, our approach achieves higher accuracy with fewer degrees of freedom, eliminates the need for complex interface meshing and fitting, and significantly reduces training time, outperforming deep neural networks.


\end{quote}

\vspace{2mm}

\section{Introduction}

Many complex physical problems, such as multiphase flow and fluid-structure interaction, require interface interactions to model multiphysics coupling (\cite{hou1997,dowell2001,gerstenberger08,khoo2009,gross11}). Therefore, the numerical treatment of the interface problem is crucial for such problems. Various numerical methods have been developed to tackle interface problems and achieve some success, such as the immersed boundary method (\cite{peskin77}), the interface-fitted finite element method (\cite{xu1982,Bramble1996,chen1998}), the immersed interface method (\cite{Li1994,li2006,lin2015}), the extended FEM (\cite{hansbo02,wang2018,xiao2020}), the discontinuous Galerkin method (\cite{Massjung2012}), the weak Galerkin method (\cite{mu2013}), and the virtual element method (\cite{chen2017}). However, these numerical methods often encounter difficulties in dealing with the coupling conditions on the interface, which may have complex geometry and dynamic boundaries. Moreover, these numerical methods face difficulties to handle high-dimensional problems. 

Neural network-based methods for solving partial differential equations (PDEs) have gained popularity in recent years, thanks to their universal approximation property (\cite{barron1993,chen1995}). Examples of such methods include the Deep Ritz Method (\cite{e2017}), the Deep Galerkin Method (\cite{sirignano2018}), the Physics-Informed Neural Networks (PINNs, \cite{raissi2019}), and so on. These methods are mesh-free and have demonstrated remarkable capabilities for solving high-dimensional PDEs and handling problems with complex geometry. However, conventional neural network training processes are inefficient because they involve solving nonlinear optimization problems that are time-consuming and may get stuck in local optima. Some deep neural networks have been proposed to solve interface problems in a mesh-free and flexible manner (\cite{wang2020,zhu2022,hu2020}), but they still face challenges in terms of training efficiency and accuracy.

Randomized Neural Networks (RNNs, \cite{pao1994,pao1995}), which do not need optimization solvers during training, are a way to overcome this challenge. Extreme Learning Machine (ELM, \cite{huang2006,huang2011}) is a special type of RNN that has been applied in various contexts (\cite{Balasundaram2011,huang2012,Sun2018,Dwivedi2020}). RNN works by randomly selecting and fixing all biases and weights except for those in the last layer, resulting in a linear combination. The weights of the last layer are then computed by solving a least-square problem. Lin et al. showed that ELM can generalize as well as fully parameterized NNs when the activation functions and initialization strategies are chosen properly (\cite{lin2014}). Based on RNN, some new methods have been proposed recently. Dong and Li developed a neural network-based method that combines local ELMs and domain decomposition to solve linear and nonlinear PDEs (\cite{dong2021}). Chen et al. extended this idea to overlapping domain decomposition (\cite{chen2022randomfm}). Shang et al. and Sun et al. incorporated RNNs with weak formulations, using the Petrov-Galerkin method (\cite{shang2022,shang2023}) and discontinuous Galerkin method (\cite{sun2022,sun2023}) respectively, to solve PDEs. These methods can achieve higher accuracy with fewer degrees of freedom and can solve time-dependent problems in the space-time approach precisely and efficiently, which suggests that this new approach has great potential for solving PDEs.

In this paper, we propose a method for solving interface problems with multiple RNNs. The whole domain is partitioned into several subdomains by the interfaces, and the solution on each subdomain is approximated by one RNN. No stochastic gradient decent type training process is required, as we obtain the solution by a least-square method, which is easier to solve than an optimization problem. This method improves the accuracy and reduces the computational cost of the numerical solution, and it can handle diffusion coefficients with large variations. Furthermore, we apply a space-time approach to solve parabolic interface problems, which avoids time steps iteration and accumulation errors.

The paper is structured as follows. Section \ref{sec2} presents the model problems. Section \ref{sec3} describes the LRNNs methods and the general process for solving interface problems. Section \ref{sec4} shows numerical results for elliptic and parabolic interface problems, illustrating the effectiveness and robustness of the proposed method. Section \ref{sec5} concludes the paper and discusses some future directions.

\section{Interface problems}
\label{sec2}

\begin{figure}[!htbp] 		
	\centering
	\includegraphics[scale=0.33]{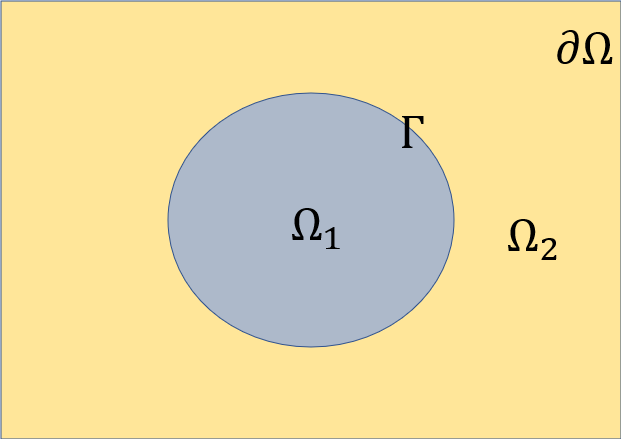}
	\caption{Domain and subdomains of an interface problem with a single closed interface.}
	\label{figure1}
\end{figure}

We present the interface problems in this section. For simplicity, we consider the case where a single closed interface in $\mathbb{R}^{2}$ splits the domain into two subdomains, as shown in Figure \ref{figure1}. Here, $\Omega\subset\mathbb{R}^{2}$ is a bounded domain with $\Omega=\Omega_1\cup\Omega_2$ and $\Gamma=\partial\Omega_1\cap\partial\Omega_2$. Our method can handle more general cases, such as multiple interfaces and higher dimensions ($\Omega \in \mathbb{R}^{d},d>2$).

\subsection{Elliptic interface problem}

We consider the following elliptic interface problem: 
\begin{equation} \label{model}
\left\{
\begin{array}{rrll}
    -\nabla \cdot (\beta(\mathbf{x}) \nabla u) &=& f, &\mathbf{x} \in \Omega, \\
    \left[u\right] &=& g_{1}, &\mathbf{x} \in \Gamma, \\
    \left[\beta(\mathbf{x}) \nabla u \cdot \mathbf{n}\right] &=& g_{2},  &\mathbf{x} \in \Gamma, \\
    u &=& g_{D},&\mathbf{x} \in \partial \Omega,
\end{array}
\right.
\end{equation}
where the diffusion coefficient $\beta(\mathbf{x})$ is a piecewise positive constant function given by
$$ \beta(\mathbf{x})=\left\{
\begin{aligned}
\beta_{1}, \quad\mathbf{x} \in \Omega_{1}, \\
\beta_{2}, \quad\mathbf{x} \in \Omega_{2}.
\end{aligned}
\right.
$$
The function $u$ is the solution of the interface problem, let $u^i$ ($i=1,2$) be the restriction of $u$ to $\Omega_i$, i.e.,
$$ 
u(\mathbf{x})=\left\{
\begin{aligned}
u^{1}(\mathbf{x}), \quad\mathbf{x} \in \Omega_{1}, \\
u^{2}(\mathbf{x}), \quad\mathbf{x} \in \Omega_{2}.
\end{aligned}
\right.
$$
The symbol $[u]$ denotes the jump of $u$ across the interface $\Gamma$, that is
$$[u]= u^1- u^2.$$
The vector $\mathbf{n}$ is the unit outward normal vector on $\Gamma$ from $\Omega_1$ to $\Omega_2$.
The function $f$ is the source term, $g_{1}$ and $g_{2}$ are the jump conditions across the interface $\Gamma$, $g_{D}$ is the boundary condition on the domain boundary $\partial\Omega$. These functions satisfy appropriate regularity conditions.

\subsection{Parabolic interface problem}

Similar to the set-up of the elliptic interface problem \eqref{model}, we consider a parabolic interface problem that involves a time interval $(0,T)$. The parabolic interface problem can be written as:
\begin{equation} \label{model_heat}
\left\{
\begin{array}{rrll}
    \frac{\partial u}{\partial t}-\nabla \cdot (\beta(\mathbf{x}) \nabla u)&=&f, \quad& (\mathbf{x},t)\in \Omega \times (0,T),\\
    \left[u\right] &=& g_{1}, \quad& (\mathbf{x},t)\in \Gamma \times (0,T),\\
    \left[\beta(\mathbf{x}) \nabla u \cdot \mathbf{n}\right] &=& g_{2},  \quad& (\mathbf{x},t)\in \Gamma \times (0,T),\\
    u(\mathbf{x},t)&=&g(\mathbf{x},t),\quad& (\mathbf{x},t)\in \partial \Omega \times (0,T),\\
        u(\mathbf{x},0)&=&u_{0}(\mathbf{x}), \quad&\mathbf{x}\in \Omega,\\
\end{array}
\right.
\end{equation}
where $u_0$ is the initial condition.
Unlike the elliptic interface problem, the interface $\Gamma$ may be time-dependent, which means that it can move over time. We assume that the interface motion is known. The problem can model a heat conduction in a composite material with a dynamic interface.

\section{Local randomized neural networks methods}
\label{sec3}

In this section, we begin by presenting the concept of randomized neural networks and their advantages. Next, we describe the LRNNs method and the mixed LRNNs method for solving elliptic interface problems. Finally, we propose the space-time LRNNs method for tackling the parabolic interface problem.

\subsection{Randomized neural networks}

We consider a fully connected neural network $\Psi: \mathbb{R}^{n_0}\rightarrow \mathbb{R}^{n_D}$ with depth $D$, which is defined as follows: 
\begin{align*}
    &\Psi_{0}(\mathbf{x}) = \mathbf{x},\\
    &\Psi_{l}(\mathbf{x}) = \rho(W_{l}\Psi_{l-1}+b_{l}), \quad l=1,...,D-1,\\
    &\Psi = \Psi_{D} = W_{D}\Psi_{D-1},
\end{align*}
where $\rho$ is the activation function, and $W_{l}\in \mathbb{R}^{n_{l}}\times \mathbb{R}^{n_{l-1}}$ and $b_{l}\in \mathbb{R}^{n_l}$ are the weights and the bias in the $l$-th layer. 
We denote the set of functions that can be represented by the above neural network by $\mathcal{N}_{\rho}^D$.
Note that $\mathcal{N}_{\rho}^D$ is not a vector space.

\begin{figure}[!htbp] 		
	\centering
	\includegraphics[scale=0.6]{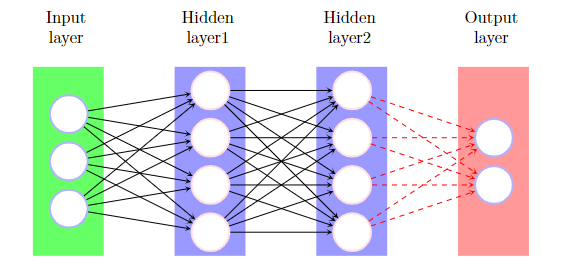}
	\caption{The structure of an RNN: the weights connecting the input layer, the first and the second hidden layers (denoted by black solid line) are randomly initialized and fixed, while the weights between the last hidden layer and output layer (denoted by red dotted line) are tunable.}
	\label{RNN}
\end{figure}

A randomized neural network (RNN) is a special type of NN, where the weights and bias in the hidden layers are randomly generated and fixed, and only the weights of the output layer are adjustable. Figure \ref{RNN} shows the structure of a two-hidden layers RNN. Note that $\mathcal{N}_{\rho}^D$ is a function space for RNNs, and $\{\psi_{1},..., \psi_{n_{D-1}} \}$ is a basis for this space, where $\psi_{i}$ is the $i$-th component of the vector $\Psi_{D-1}$. Hence, any function $u_{\rho}\in\mathcal{N}_{\rho}^D$ (one neuron in the output layer) can be expressed as a linear combination of the basis functions
\begin{equation}
    u_{\rho} = \sum\limits_{i=1}^{n_{D-1}} \alpha_{i} \psi_{i}. \label{linear_combination}
\end{equation}

\subsection{Local RNNs method for elliptic interface problem}

The domain $\Omega$ is partitioned into several subdomains by interfaces, and the solutions in different subdomains may have distinct characteristics, so one neural network may not be adequate to describe the solution in the whole domain. Therefore, we use multiple RNNs to learn interface problems. Specifically, we use one RNN for each subdomain, i.e., the number of RNNs equals the number of subdomains.

Consider the case where the domain $\Omega$ is split into two subdomains by interface $\Gamma$. We propose a local RNNs method to solve the elliptic interface problem \eqref{model}. In this method, we use two RNNs to approximate the true solution in each subdomain denoted by $u^{1}_{\rho}$ and $u^{2}_{\rho}$, which correspond to the subdomains $\Omega_{1}$ and $\Omega_{2}$. For simplicity, we assume that both RNNs have the same number of neurons $m$ in the hidden layer. Similar to \eqref{linear_combination}, we assume that the outputs of the RNNs are given by
\begin{align}
    &u^{1}_{\rho} = \sum\limits_{k=1}^{m} \alpha^{1}_{k} \psi^{1}_{k}, \qquad
    u^{2}_{\rho} = \sum\limits_{k=1}^{m} \alpha^{2}_{k} \psi^{2}_{k}. \label{combination11}
\end{align}
If we substitute \eqref{combination11} into the elliptic interface problem \eqref{model}, and discretize the problem at randomly sampled $N_1$, $N_2$ and $N_3$ collocation points on $\Omega$, $\Gamma$ and $\partial \Omega$, respectively, we obtain a linear system
\begin{equation}\label{boundary_equation}
\begin{array}{ll}
    &A\mathit{X} = F,\\
    &B\mathit{X} = G^{1},\\
    &C\mathit{X} = G^{2},\\
    &D\mathit{X} = G^{D},
\end{array}
\end{equation}
where $A,B,C,D$ are matrices of order $N_{1}\times 2m$, $N_{2}\times 2m$, $N_{2}\times 2m$ and $N_{3}\times 2m$, respectively. Specifically, 
$$ A_{i,j}=\left\{
\begin{aligned}
-\nabla \cdot (\beta(\mathbf{x}_{i}) \nabla \psi^{1}_{j}(\mathbf{x}_{i})), \quad &\mathbf{x}_{i} \in \Omega_{1},\; 1\le j\le m, \\
-\nabla \cdot (\beta(\mathbf{x}_{i}) \nabla \psi^{2}_{j-m}(\mathbf{x}_{i})), \quad &\mathbf{x}_{i} \in \Omega_{2},\; m+1 \le j\le 2m, 
\end{aligned}
\right.
$$
$$
B_{i,j}=\left\{
\begin{aligned}
\psi^{1}_{j}(\mathbf{x}_{i}), \quad &\mathbf{x}_{i} \in \Gamma,\;  1\le j\le m, \\
-\psi^{2}_{j-m}(\mathbf{x}_{i}), \quad &\mathbf{x}_{i} \in \Gamma,\;  m+1 \le j\le 2m,
\end{aligned}
\right.
$$
$$ C_{i,j}=\left\{
\begin{aligned}
\nabla \psi^1_{j}\cdot \mathbf{n}(\mathbf{x}_{i}), \quad &\mathbf{x}_{i} \in \Gamma,\;  1\le j\le m, \\
-\nabla \psi^2_{j-m}\cdot \mathbf{n}(\mathbf{x}_{i}), \quad &\mathbf{x}_{i} \in \Gamma,\;  m+1 \le j\le 2m,
\end{aligned}
\right.
\quad D_{i,j}=\left\{
\begin{aligned}
0, \quad &\mathbf{x}_{i} \in \partial\Omega,\; 1\le j\le m, \\
\psi^{2}_{j-m}(\mathbf{x}_{i}), \quad &\mathbf{x}_{i} \in \partial\Omega,\;  m+1 \le j\le 2m,
\end{aligned}
\right.
$$
$$
F_i = f(\mathbf{x}_{i}), \quad G^1_i = g_1(\mathbf{x}_{i}),\quad G^2_i = g_2(\mathbf{x}_{i}), \quad G_i^D = g_D(\mathbf{x}_{i}),
$$
and
$\mathit{X}=(\alpha^{1}_{1},...,\alpha^{1}_{m},\alpha^{2}_{1},...,\alpha^{2}_{m})^{T}$ is the unknown vector. 

Automatic differentiation is a method to compute derivatives of functions, but it can be inefficient for some applications. Therefore, we use finite difference schemes to approximate derivatives. For instance, we can use the 1st and 2nd order central difference schemes to approximate the 1st and 2nd order partial derivatives of $u_\rho(x,y)$ with respect to x, as shown below
\begin{align}
	\frac{\partial{u_\rho(x,y)}}{\partial x} &= \frac{u_\rho(x+h_x,y)-u_\rho(x-h_x,y)}{2h_x}, \label{diff_1}\\
	\frac{\partial^2{u_\rho(x,y)}}{\partial^2 x} &= \frac{u_\rho(x+h_x,y) -2u_\rho(x,y)  + u_\rho(x-h_x,y) }{h_x^2}. \label{diff_2}
\end{align}
Here, $h_x$ is the step size.
By solving a least-square problem with the linear system \eqref{boundary_equation}, we can obtain the weights of the output layer for the RNNs.

\begin{remark}
The local RNNs method has some advantages over the traditional NN method because we can train the RNNs by a least-squares method, which avoids the issues of local minima and time-consuming training process that may arise from solving a nonlinear optimization problem.
\end{remark}

\begin{remark}
We can adjust the weights $\gamma_{1}$, $\gamma_{2}$ and $\gamma_{3}$ in the linear system to account for the importance of different conditions:
\begin{align*}
    &A\mathit{X} = F,\\
    &\gamma_{1} B\mathit{X} = \gamma_{1} G^{\Gamma1},\\
    &\gamma_{2} C\mathit{X} = \gamma_{2} G^{\Gamma2},\\
    &\gamma_{3} D\mathit{X} = \gamma_{3} G^{D}.
\end{align*}
We can choose different weights to control the attention of RNNs. For example, we can assign larger weights to the interface conditions and boundary conditions, which are usually more crucial for the accuracy of the solution.
\end{remark}

\subsection{Mixed LRNNs method for elliptic interface problem}

By introducing another vector-valued function $\mathbf{p} = \beta(\mathbf{x}) \nabla u$, we rewrite problem \eqref{model}
as following mixed form:
\begin{equation} \label{mix}
\left\{
\begin{array}{rrll}
  -\nabla \cdot \mathbf{p} &=& f,\quad &\mathbf{x} \in \Omega, \\
  \mathbf{p} - \beta(\mathbf{x}) \nabla u &=& 0,\quad&\mathbf{x} \in \Omega,\\
  \left[u\right] &=& g_{1},\quad&\mathbf{x} \in \Gamma, \\
  \left[\mathbf{p} \cdot \mathbf{n}\right] &=& g_{2},\quad&\mathbf{x} \in \Gamma, \\
  u &=& g_{D},\quad&\mathbf{x} \in \partial \Omega,
\end{array}
\right.
\end{equation}
where all notations are the same as before. More RNNs are needed to approximate $\mathbf{p}$. Therefore, we set their outputs as follows:
\begin{align*}
  &u^{1}_{\rho} = \sum\limits_{k=1}^{m} \alpha^{1}_{k} \psi^{1}_{k},\qquad
  u^{2}_{\rho} = \sum\limits_{k=1}^{m} \alpha^{2}_{k} \psi^{2}_{k},\\
  &\mathbf{p}^{1}_{\rho} =\left( \sum\limits_{k=1}^{m} \tau^{1,1}_{k} \xi^{1,1}_{k},\; \sum\limits_{k=1}^{m} \tau^{1,2}_{k} \xi^{1,2}_{k}\right),\qquad
  \mathbf{p}^{2}_{\rho} =\left( \sum\limits_{k=1}^{m} \tau^{2,1}_{k} \xi^{2,1}_{k},\; \sum\limits_{k=1}^{m} \tau^{2,2}_{k} \xi^{2,2}_{k}\right).
\end{align*}
Substitute them into the mixed elliptic interface problem \eqref{mix}, and add weights $\gamma_{i}$ ($i=1,2,3$), we achieve the linear system: 
\begin{align*}
  &A^{1}\mathit{X} = F,\\
  &A^{2}\mathit{X} = 0,\\
  &A^{3}\mathit{X} = 0,\\
  &\gamma_{1} B\mathit{X} = \gamma_{1} G^{1},\\
  &\gamma_{2} C\mathit{X} = \gamma_{2} G^{2},\\
  &\gamma_{3} D\mathit{X} = \gamma_{3} G^{D},
\end{align*}
where $A^{1},A^{2},A^{3},B,C,D$ are matrices of order $N_{1}\times 6m, N_{1}\times 6m, N_{1}\times 6m, N_{2}\times 6m, N_{2}\times 6m, N_{3}\times 6m$, respectively. Here,

$$ A^{1}_{i,j}=\left\{
\begin{aligned}
-\frac{\partial \xi^{1,1}_{j-2m}}{\partial x}(\mathbf{x}_{i}), \quad \mathbf{x}_{i} \in \Omega_{1},\; 2m+1\le j\le 3m, \\
-\frac{\partial \xi^{1,2}_{j-3m}}{\partial y}(\mathbf{x}_{i}), \quad \mathbf{x}_{i} \in \Omega_{1},\; 3m+1\le j\le 4m, \\
-\frac{\partial \xi^{2,1}_{j-4m}}{\partial x}(\mathbf{x}_{i}), \quad \mathbf{x}_{i} \in \Omega_{2},\; 4m+1\le j\le 5m, \\
-\frac{\partial \xi^{2,2}_{j-5m}}{\partial y}(\mathbf{x}_{i}), \quad \mathbf{x}_{i} \in \Omega_{2},\; 5m+1\le j\le 6m,
\end{aligned}
\right.
\quad
 A^{2}_{i,j}=\left\{
\begin{aligned}
-\beta_{1}\frac{\partial \psi^{1}_{j}}{\partial x}(\mathbf{x}_{i}), \quad & \mathbf{x} \in \Omega_{1},\; 1\le j\le m, \\
-\beta_{2}\frac{\partial \psi^{2}_{j-m}}{\partial x}(\mathbf{x}_{i}), \quad &\mathbf{x} \in \Omega_{2},\; m+1\le j\le 2m, \\
\xi^{1,1}_{j-2m}(\mathbf{x}_{i}), \quad &\mathbf{x} \in \Omega_{1},\; 2m+1\le j\le 3m, \\
\xi^{2,1}_{j-4m}(\mathbf{x}_{i}), \quad &\mathbf{x} \in \Omega_{2},\; 4m+1\le j\le 5m, \\
\end{aligned}
\right.
$$
$$ A^{3}_{i,j}=\left\{
\begin{aligned}
-\beta_{1}\frac{\partial \psi^{1}_{j}}{\partial y}(\mathbf{x}_{i}), \quad &\mathbf{x}_{i} \in \Omega_{1},\; 1\le j\le m, \\
-\beta_{2}\frac{\partial \psi^{2}_{j-m}}{\partial y}(\mathbf{x}_{i}), \quad &\mathbf{x}_{i} \in \Omega_{2},\; m+1\le j\le 2m, \\
\xi^{1,2}_{j-3m}(\mathbf{x}_{i}), \quad &\mathbf{x}_{i} \in \Omega_{1},\; 3m+1\le j\le 4m, \\
\xi^{2,2}_{j-5m}(\mathbf{x}_{i}), \quad &\mathbf{x}_{i} \in \Omega_{2},\; 5m+1\le j\le 6m, \\
\end{aligned}
\right.
\quad
C_{i,j}=\left\{
\begin{aligned}
(\xi^{1,1}_{j-2m}n_{1})(\mathbf{x}_{i}), \quad\mathbf{x}_{i} \in \Gamma,\;  2m+1\le j\le 3m, \\
(\xi^{1,2}_{j-3m}n_{2})(\mathbf{x}_{i}), \quad\mathbf{x}_{i} \in \Gamma,\;  3m+1\le j\le 4m,\\
-(\xi^{2,1}_{j-4m}n_{1})(\mathbf{x}_{i}), \quad\mathbf{x}_{i} \in \Gamma,\;  4m+1\le j\le 5m, \\
-(\xi^{2,2}_{j-5m}n_{2})(\mathbf{x}_{i}), \quad\mathbf{x}_{i} \in \Gamma,\;  5m+1\le j\le 6m,
\end{aligned}
\right.
$$
$$ B_{i,j}=\left\{
\begin{aligned}
\psi^{1}_{j}(\mathbf{x}_{i}), \quad  & \mathbf{x}_{i} \in \Gamma,\; 1\le j\le m, \\
-\psi^{2}_{j-m}(\mathbf{x}_{i}), \quad & \mathbf{x}_{i} \in \Gamma,\; m+1 \le j\le 2m,
\end{aligned}
\right.
\quad
 D_{i,j}=\left\{
\begin{aligned}
\psi^{2}_{j-m}(\mathbf{x}_{i}), \quad & \mathbf{x}_{i} \in \partial\Omega,\;   m+1 \le j\le 2m,\\
0, \quad &\mathbf{x}_{i} \in \partial\Omega,\;   {\rm otherwise},
\end{aligned}
\right.
$$
and $\mathit{X}=(\alpha^{1}_{1},\cdots,\alpha^{1}_{m},\alpha^{2}_{1},...,\alpha^{2}_{m},\tau^{1,1}_{1},...,\tau^{1,1}_{m},\tau^{1,2}_{1},\cdots,\tau^{1,2}_{m},\tau^{2,1}_{1},...,\tau^{2,1}_{m},\tau^{2,2}_{1},...,\tau^{2,2}_{m})$.

\subsection{Space-time LRNNs method for parabolic interface problem}

We use a space-time approach to solve the parabolic interface problem \eqref{model_heat}. In the space-time approach, the time variable and space variables are treated equally. We assume that the interface has a smooth and known motion. When $\Omega\subset \mathbb{R}^2$, the interface becomes a surface in the 3-dimensional space-time domain, dividing it into two subdomains. The initial condition can be viewed as a boundary condition for the space-time domain. We employ two RNNs to approximate the solution in each sub-domain.

\begin{align}
  &u^{1}_{\rho} = \sum\limits_{k=1}^{m} \alpha^{1}_{k} \psi^{1}_{k}, \qquad
  u^{2}_{\rho} = \sum\limits_{k=1}^{m} \alpha^{2}_{k} \psi^{2}_{k}. \label{combination22}
\end{align}
Then we have:
\begin{align}\label{linear system}
  &A\mathit{X} = F,\\
  &\gamma_{1} B\mathit{X} = \gamma_{1} G^{\Gamma1},\\
  &\gamma_{2} C\mathit{X} = \gamma_{2} G^{\Gamma2},\\
  &\gamma_{3} D\mathit{X} = \gamma_{4} G^{D},\\
  &\gamma_{4} E\mathit{X} = \gamma_{3} U^{0},\\
\end{align}
where $A,B,C,D,E$ are $N_{1}\times 2m, N_{2}\times 2m, N_{2}\times 2m, \frac{4N_{3}}{5}\times 2m, \frac{N_{3}}{5}\times 2m$-order matrix, respectively,

$$ A_{i,j}=\left\{
\begin{aligned}
\left(\frac{\partial \psi^{1}_{j}}{t}-\nabla \cdot (\beta(\mathbf{x}_{i}) \nabla\psi^{1}_{j})\right)(\mathbf{x}_i,t),\quad &\mathbf{x}_i \in \Omega_{1},\; 1\le j\le m, \\
\left(\frac{\partial \psi^{2}_{j-m}}{t}-\nabla \cdot (\beta(\mathbf{x}_{i}) \nabla\psi^{2}_{j -m})\right)(\mathbf{x}_i,t),\quad &\mathbf{x}_i \in \Omega_{2},\; m+1 \le j\le 2m, 
\end{aligned}
\right.
$$
$$ B_{i,j}=\left\{
\begin{aligned}
\psi^{1}_{j}(\mathbf{x}_i,t),\quad &\mathbf{x}_{i} \in \Gamma,\;   1\le j\le m, \\
-\psi^{2}_{j-m}(\mathbf{x}_i,t),\quad & \mathbf{x}_{i} \in \Gamma,\;   m+1 \le j\le 2m,
\end{aligned}
\right.
\quad
 C_{i,j}=\left\{
\begin{aligned}
(\nabla \psi^1_{j}\cdot \mathbf{n})(\mathbf{x}_i,t),\quad & \mathbf{x}_{i} \in \Gamma,\;   1\le j\le m, \\
(-\nabla \psi^2_{j-m}\cdot \mathbf{n})(\mathbf{x}_i,t),\quad  & \mathbf{x}_{i} \in \Gamma,\;   m+1 \le j\le 2m,
\end{aligned}
\right.
$$
$$ D_{i,j}=\left\{
\begin{aligned}
0,\quad & \mathbf{x}_{i} \in \partial\Omega,\;   1\le j\le m, \\
\psi^{2}_{j-m}(\mathbf{x}_i,t),\quad & \mathbf{x}_{i} \in \partial\Omega,\;   m+1 \le j\le 2m,
\end{aligned}
\right.
\quad
 E_{i,j}=\left\{
\begin{aligned}
\psi^{1}_{j}(\mathbf{x}_i,0),\quad & \mathbf{x}_i \in \Omega_{1},\; 1\le j\le m, \\
\psi^{2}_{j-m}(\mathbf{x}_i,0),\quad & \mathbf{x}_i \in \Omega_{2},\; m+1\le j\le 2m,
\end{aligned}
\right.
$$
and $X=(\alpha^{1}_{1},...,\alpha^{1}_{m},\alpha^{2}_{1},...,\alpha^{2}_{m})^{T}$.

\begin{remark}
Breaking down a large problem into smaller pieces and solving them separately can be a good way to save memory and time for many problems. That is why finite difference methods are popular for the time variable, which follow this strategy. But this method has a problem: it can make errors that get bigger and bigger over time. The LRNNs method does not have this problem because it uses a space-time approach. It can solve time-dependent problems efficiently and accurately by solving one least square problem.
\end{remark}

\section{Numerical examples}
\label{sec4}

We demonstrate the performance of the local randomized neural network methods on several test problems in this section.

We approximate solutions using single-hidden layer RNNs with $\rho=\tanh$ as the activation function and finite difference schemes to approximate the gradient and Laplacian operators. The step sizes are $h_x = 10^{-6}$ and $h_x = 5\times 10^{-4}$ for the 1st and 2nd order partial derivatives, respectively, in \eqref{diff_1} and \eqref{diff_2}. We draw the initial weights and bias from the uniform distribution $\mathcal{U}(-r_i, r_i)$ for each layer and fix them except for the last layer, where $r\in \mathbb{R}$. We compute the $L^{2}(\Omega)$  errors in Example \ref{ex1} – \ref{ex5} using $20^{d}$ Gauss-Legendre quadrature points, where $d$ is the dimension of the domain (or space-time domain). We set $\gamma_i=50$ for all weights in the linear system. We average the results over ten experiments for each problem to account for the randomness in the methods. We use PyTorch to implement the numerical experiments.

\begin{example}\label{ex1}
    {Two-dimensional Flower Shape Interface Problem.}
Let $\Omega = (-1,1)^2$, we consider an elliptic interface problem \eqref{model} with a flower shape interface
$$
\Gamma:(x-x_{c})^2+(y-y_{c})^2=r^{2}(\theta),\quad r(\theta)=0.4+0.2\sin(20\theta),\quad \theta \in [0,2\pi),
$$
where $x_{c}=y_{c}=0.02\sqrt{5}$. The exact solution is given by
$$ u(\mathbf{x})=\left\{
\begin{aligned}
&\frac{1}{\beta_{1}}e^{xy},\quad \mathbf{x}\in \Omega_{1}, \\
&\frac{1}{\beta_{2}}\sin(x)\sin(y),\quad \mathbf{x}\in \Omega_{2},
\end{aligned}
\right.
$$
$f$, $g_{1}$ and $g_{2}$ can be calculated from the exact solution and coefficients.
\end{example}

To learn the solution, we use two RNNs with the same structure: the input layer has $n_0 = 2$ neurons, the output layer has  $n_2 = 1$ neuron and the hidden layer has $n_1 = m$ neurons. The weights and bias are initialized from the uniform distributions $\mathcal{U}(-r_{1},r_{1})$ and $\mathcal{U}(-r_{2},r_{2})$, respectively. We may choose different ranges of $r_1$ and $r_2$ for different cases.

\begin{figure}[!htbp]
      \centering     
      \subfloat[approximation$u_\rho$]
      {
        \label{flower_1}\includegraphics[width=0.33\textwidth]{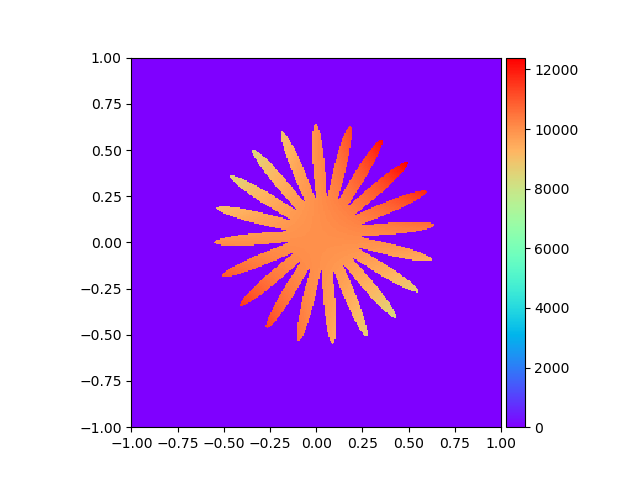}
      }
      \subfloat[real solution$u$]
      {
        \label{flower_2}\includegraphics[width=0.33\textwidth]{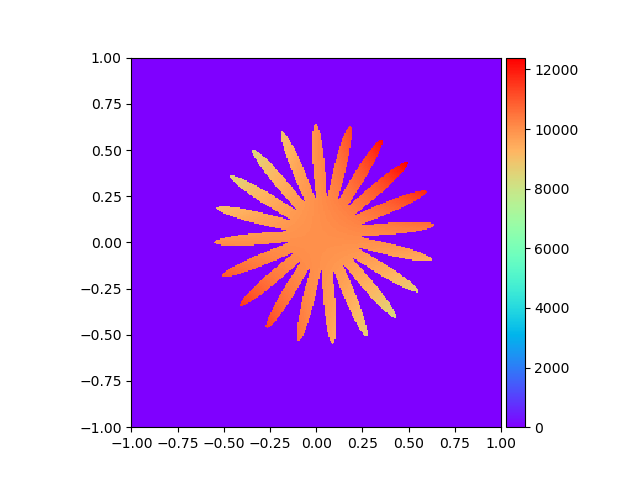}
      }
      \subfloat[absolute error$|u-u_{\rho}|$]
      {
        \label{flower_3}\includegraphics[width=0.33\textwidth]{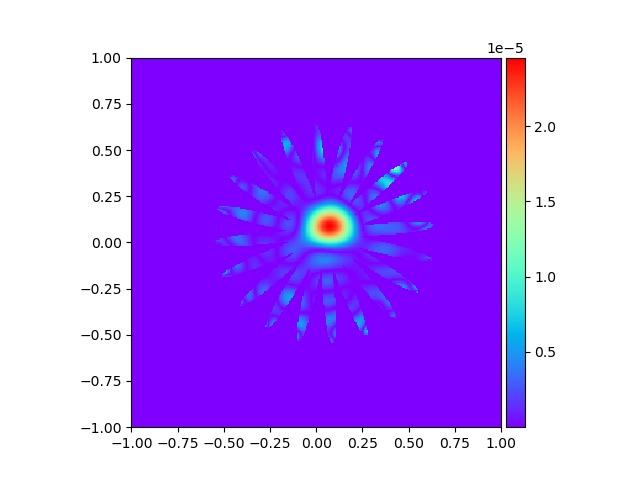}
      }
      \caption{Displays of $u_\rho$, $u$ and $|u-u_{\rho}|$ for LRNNs method with $\beta_{1}=0.0001$, $\beta_{2}=10000$, $r_{1}=1.6$, $r_{2}=0.7$, $m=320$, $N=5000$ in Example \ref{ex1}}.
      \label{figure2}

\end{figure}

Figure \ref{figure2} shows the numerical solution $u_\rho$, the exact solution $u$ and their error $|u-u_{\rho}|$.  In this case, we choose the coefficients $\beta_{1}=0.0001$ and $\beta_{2}=10000$, set $m=320$, $r_{1}=1.6$, $r_{2}=0.7$, and sample $5000$ points uniformly from $\Omega$, $\Gamma$ and $\partial \Omega$ with the ratio of $3:1:1$. 

Table \ref{tab1} shows the relative $L^2$ errors for the LRNNs method with different numbers of sampled points $N$ and neurons $m$ in the hidden layers. The number of degrees of freedom is $2m$, i.e., ${\rm dof}=2m$. The results reveal that the error decreases as the number of neurons increases, which enhances the expressiveness of the space $N_\rho^D$. The accuracy also improves with more sampled points, but the improvement is marginal for small values of $m$. This implies that the expressiveness of LRNNs and the data size should be balanced for optimal approximation.

\begin{table}[!htbp]	
	\centering  
	\begin{tabular}{|c|c|c|c|c|c|c|}  
		\hline  
		\diagbox [width=6em] {N}{m}&40&80&160&320&640&1280\cr\cline{1-7}
        500&1.89E-02&1.73E-03&1.94E-05&4.17E-07&2.19E-07&4.68E-07\cr\hline
        1000&1.97E-02&1.18E-03&9.20E-06&3.94E-08&2.51E-08&3.23E-08\cr\hline
		    2000&2.26E-02&1.46E-03&8.32E-06&3.15E-08&2.77E-08&1.20E-08\cr\hline
		    5000&1.92E-02&1.10E-03&9.04E-06&3.17E-08&2.51E-08&1.77E-08\cr\hline 
		    10000&2.27E-02&1.33E-03&9.17E-06&2.93E-08&2.59E-08&2.43E-08\cr\hline
	\end{tabular}  
	\caption{Relative $L^2(\Omega)$ errors of the LRNNs method with different $N$ and $m$, where $r_{1}=1.6$, $r_{2}=0.7$, $\beta_{1}=1$ and $\beta_{2}=10$ in Example \ref{ex1}.}\label{tab1}
\end{table}

The LRNNs method’s performance with varying values of $\beta_{1}$ and $\beta_{2}$ is shown in Table \ref{tab2}, where $N=5000$ and $m=320$. The results demonstrate that the LRNNs method can achieve high accuracy for interface problems with large coefficient variations. The method is also computationally efficient and has a competitive running time.

\begin{table}[!htbp]
    \centering
    \begin{tabular}{|c|c|c|c|c|c|}
    \hline
    $(\beta_{1},\beta_{2})$ & $(r_{1},r_{2})$ &relative $L^{2}$ error & CPU time (s) \\\hline
    ($1,10^{2}$)&(1.2,1.2)&1.36E-08&1.9 \\\hline
    ($1,10^{4}$)&(1.8,1.3)&1.03E-08&1.9  \\\hline
    ($10^{-4},10^{4}$)&(1.6,0.7)&1.17E-09&1.9  \\\hline
    ($10^{-6},10^{6}$)&(1,1.3)&1.20E-07&1.9  \\\hline
    ($10^{2},10^{-2}$)&(1.3,1.6)&4.05E-08&2.4  \\\hline
    ($10^{4},10^{-4}$)&(1,1.8)&6.26E-05&2.1 \\\hline
    \end{tabular}
    \caption{Relative $L^2(\Omega)$ errors of the LRNNs method for different $\beta_{1}$ and $\beta_{2}$, where $m=320$ and $N=5000$ in Example \ref{ex1}.}\label{tab2}
\end{table}

%We approximate the exact solution with four RNNs: $u_\rho^1$ and $u_\rho^2$ for $u$ in $\Omega_{1}$ and $\Omega_{2}$, respectively, and $\mathbf{p}_\rho^1$ and $\mathbf{p}_\rho^2$ for $\mathbf{p}$ in $\Omega_{1}$ and $\Omega_{2}$, respectively. We set $\beta_{1}=0.0001$, $\beta_{2}=10000$, $N=5000$, $m=320$, and sample points on $\Omega$, $\Gamma$, $\partial \Omega$ in a ratio of $3:1:1$. Figure \ref{figure10} shows the approximated solution $u_{\rho}$, the exact solution $u$, and the absolute error $|u-u_{\rho}|$.

%\begin{figure}[!htbp]
%  \centering     
%  \subfloat[approximation$u_\rho$]
%  {
%    \label{mix_1}\includegraphics[width=0.3\textwidth]{mix1.png}
%  }
%  \subfloat[real solution$u$]
%  {
%    \label{mix_2}\includegraphics[width=0.3\textwidth]{mix2.png}
%  }
%  \subfloat[absolute error$|u-u_{\rho}|$]
%  {
%    \label{mix_3}\includegraphics[width=0.3\textwidth]{mix3.png}
%  }
%  \caption{Displays of $u_\rho$, $u$ and $|u-u_{\rho}|$ for mixed LRNNs method with $\beta_{1}=0.0001$, $\beta_{2}=10000$, $r_{1}=1$, $r_{2}=1.1$, $r_{3}=0.7$, $r_{4}=2.1$, $m=320$, $N=5000$ in Example \ref{ex1}.}
%  \label{figure10}
%\end{figure}

Let us apply the mixed LRNNs method to this problem with the mixed form \eqref{mix}. 
Table \ref{tab9} shows the relative $L^2$ errors for different values of $m$ and $N$, with $\beta_{1}=1$, $\beta_{2}=10$, $r_{1}=1$, $r_{2}=1.1$, $r_{3}=0.7$ and $r_{4}=2.1$. We can see that the mixed LRNNs method has a similar performance to the LRNNs method, meaning that both methods achieve higher accuracy with more neurons and larger $N$. Moreover, we notice that the mixed LRNNs method can obtain smaller errors when the number of sampling points is sufficiently large. This is because the mixed LRNNs method uses four RNNs to approximate the same problem, so it can capture the solution better with more data points.

\begin{table}[!htbp]	
	\centering  
	\begin{tabular}{|c|c|c|c|c|c|c|c|c|}  
		\hline  
		\diagbox [width=4em] {N}{m}&  
		\multicolumn{2}{c|}{80}&\multicolumn{2}{c|}{160}&
		\multicolumn{2}{c|}{320}&\multicolumn{2}{c|}{640}\cr\cline{1-9}  
		&$\frac{\Vert u-u_{\rho}\Vert}{\Vert u\Vert_2} $ & $\frac{\Vert \mathbf{p} - \mathbf{p}_{\rho}\Vert_2}{\Vert \mathbf{p}\Vert_2}$ &$\frac{\Vert u-u_{\rho}\Vert}{\Vert u\Vert_2} $ & $\frac{\Vert \mathbf{p} - \mathbf{p}_{\rho}\Vert_2}{\Vert \mathbf{p}\Vert_2}$&$\frac{\Vert u-u_{\rho}\Vert}{\Vert u\Vert_2} $ & $\frac{\Vert \mathbf{p} - \mathbf{p}_{\rho}\Vert_2}{\Vert \mathbf{p}\Vert_2}$&$\frac{\Vert u-u_{\rho}\Vert}{\Vert u\Vert_2} $ & $\frac{\Vert \mathbf{p} - \mathbf{p}_{\rho}\Vert_2}{\Vert \mathbf{p}\Vert_2}$\cr\hline  
		500&3.02E-04&1.02E-02&1.65E-03&6.40E-03&1.27E-04&3.15E-05&1.11E-05&1.90E-05\cr\hline  
		1000&1.18E-04&2.48E-03&1.51E-06&3.36E-05&6.54E-07&1.01E-06&1.54E-07&3.22E-07\cr\hline  
		2000&2.29E-04&2.16E-03&1.08E-06&1.22E-05&1.96E-08&1.71E-07&3.75E-09&4.65E-08\cr\hline  
		5000&1.29E-04&2.02E-03&5.15E-07&9.81E-06&9.57E-10&4.04E-08&4.24E-10&8.59E-09\cr\hline  
		10000&1.66E-04&2.09E-03&4.71E-07&8.82E-06&8.98E-10&3.66E-08&1.71E-10&6.59E-09\cr\hline  
	\end{tabular}  
	\caption{Relative $L^2(\Omega)$ errors of the mixed LRNNs method with different $N$ and $m$, where $r_{1}=1$, $r_{2}=1.1$, $r_{3}=0.7$, $r_{4}=2.1$, $\beta_{1}=1$ and $\beta_{2}=10$ in Example \ref{ex1}.}\label{tab9}
\end{table}

Table \ref{tab10} shows the relative $L^2$ error and CPU time for different diffusion coefficients, with $m=320$ and $N=5000$. Unlike the results of the LRNNs method in Table \ref{tab2}, the mixed LRNNs method takes more CPU time because it has more dof. The mixed LRNNs method is also robust to the diffusion coefficients.

\begin{table}[!htbp]
  \centering
  \begin{tabular}{|c|c|c|c|c|c|}
  \hline
  $(\beta_{1},\beta_{2})$ &$(r_{1},r_{2},r_{3},r_{4})$& relative $L^{2}$ error & CPU time (s) \\\hline
  ($1,10^{2}$)&(1.2,0.8,2.6,0.9)&1.32E-09&8.7 \\\hline
  ($1,10^{4}$)&(1.9,0.6,1,1.2)&1.71E-07&8.5  \\\hline
  ($10^{-4},10^{4}$)&(1,1.1,0.7,2.1)&2.15E-10&8.3  \\\hline
  ($10^{-6},10^{6}$)&(1,0.5,0.5,1.4)&7.57E-09&8.5  \\\hline
  ($10^{2},10^{-2}$)&(0.2,0.8,1.3,1.3)&7.32E-08&8.3  \\\hline
  ($10^{4},10^{-4}$)&(2.3,0.7,1.4,1.6)&1.14E-04&8.2 \\\hline
  \end{tabular}
  \caption{Relative $L^2(\Omega)$ errors of $u_\rho$ by the mixed LRNNs method for different $\beta_{1}$ and $\beta_{2}$, where $m=320$ and $N=5000$ in Example \ref{ex1}.}\label{tab10}
\end{table}

\begin{example}\label{ex2}
  Three-dimensional Interface Problem. Let $\Omega = (-1,1)^3$,
we consider an elliptic interface problem \eqref{model} with spherical interface
$$
\Gamma:x^2+y^2+z^{2}=0.75^{2}.
$$
The exact solution
$$ 
u(\mathbf{x})=\left\{
\begin{aligned}
5e^{x^{2}+y^{2}+z^{2}}+20,\quad\mathbf{x}\in \Omega_{1}, \\
10(x+y+z),\quad\mathbf{x}\in \Omega_{2},
\end{aligned}
\right.
$$
$f$, $g_{1}$ and $g_{2}$ can be derived from the exact solution $u$ and coefficients.
\end{example}

Since $\Gamma$ divides $\Omega$ into two subdomains, we use two RNNs with the same structure to approximate the solution. The RNNs have two layers with $n_0 = 3$, $n_1 = m$ and $n_2 = 1$. The weights and biases are randomly chosen from $\mathcal{U}(-2.54,2.54)$ and $\mathcal{U}(-0.33,0.33)$, respectively. We sample $N$ points from $\Omega$, $\Gamma$ and $\partial \Omega$ with a proportion of 6:1:3.

\begin{figure}[!htbp]
      \centering
      \subfloat[approximation $u_\rho$]
      {
        \label{sphere_1}\includegraphics[width=0.33\textwidth]{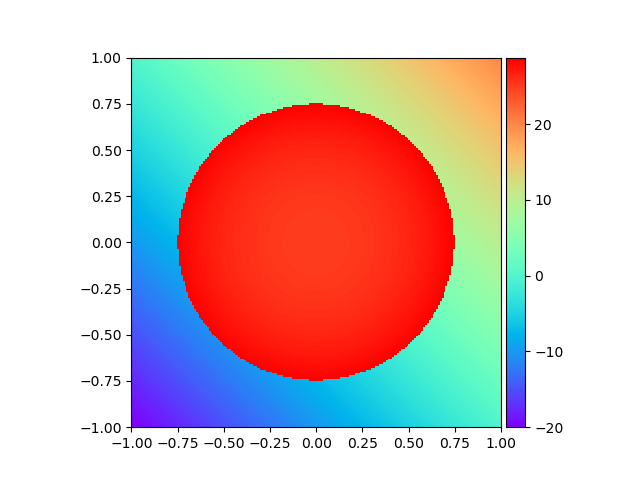}
      }
      \subfloat[real solution $u$]
      {
        \label{sphere_2}\includegraphics[width=0.33\textwidth]{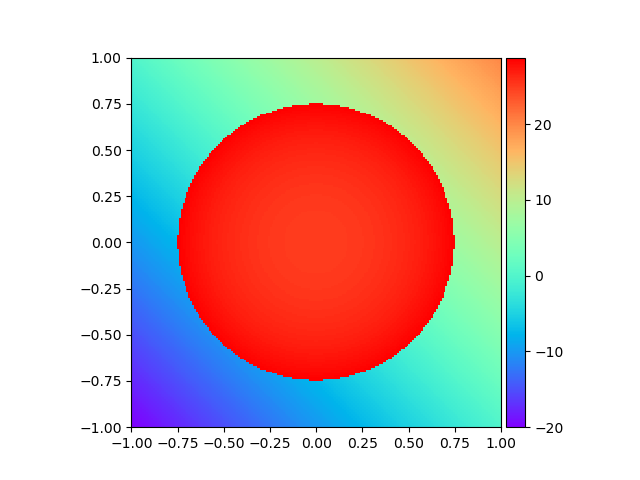}
      }
      \subfloat[absolute error $|u-u_{\rho}|$]
      {
        \label{sphere_3}\includegraphics[width=0.33\textwidth]{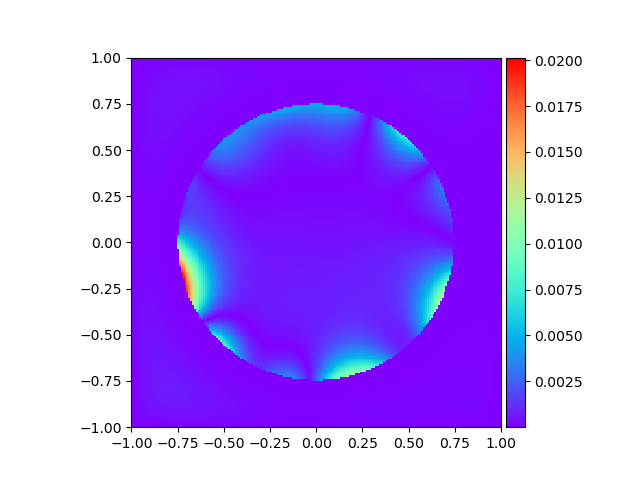}
      }
      \caption{Displays of $u_\rho$, $u$ and $|u-u_{\rho}|$ at $z=0$ for the LRNNs method with $\beta_{1}=1$, $\beta_{2}=1000$, $m=640$, $N=10000$, $r_{1}=2.54$, $r_{2}=0.33$ in Example \ref{ex2}.}
      \label{figure3}
\end{figure}

Figure \ref{figure3} shows the RNN solution $u_{\rho}$, the exact solution $u$, and the absolute error $|u-u_{\rho}|$ at $z=0$ for $\beta_{1}=1$ and $\beta_{2}=1000$. We use $m=640$, $N=10000$.
Table \ref{tab3} reports the relative $L^2$ errors for different values of $m$ and $N$. As in Example \ref{ex1}, 
 increasing $m$ and $N$ improves the numerical solution.

\begin{table}[!htbp]	
	\centering  
	\begin{tabular}{|c|c|c|c|c|c|}  
		\hline  
		\diagbox [width=6em] {N}{m}&80&160&320&640&1280\cr\cline{1-6}
        1000&4.04E-03&1.66E-03&1.10E-03&1.22E-03&7.95E-04\cr\hline
		    2000&4.31E-03&6.86E-04&3.90E-04&2.82E-04&1.25E-04\cr\hline
		    5000&4.03E-03&5.63E-04&2.19E-04&9.87E-05&1.67E-05\cr\hline 
		    10000&3.86E-03&6.66E-04&1.27E-04&3.13E-05&1.93E-05\cr\hline
        20000&3.67E-03&6.61E-04&1.46E-04&1.84E-05&1.19E-06\cr\hline
	\end{tabular}  
	\caption{Relative $L^2(\Omega)$ errors with different $N$ and $m$, where $\beta_{1}=1$, $\beta_{2}=100$ in Example \ref{ex2}.}
	\label{tab3}
\end{table}

Table \ref{tab4} shows the relative $L^2$ errors and CPU time for different diffusion coefficients $\beta_{1}$, $\beta_{2}$. We observe that the LRNNs method is robust to the diffusion coefficients. Consider the case where $\beta_{1}=1$, $\beta_{2}=10$, and compare our method with the virtual element method (Example 5.1 in \cite{chen2017}). We obtain a $L^2$ error of $2.37\times 10^{-3}$, which is comparable to their error. However, our method takes only about 6.8 seconds, while their method requires 22.4 seconds, even with the assistance of an algebraic multigrid solver for acceleration. This result highlights the superior efficiency of our approach.

%Consider the case where $\beta_{1}=1,\beta_{2}=10$, we give relative $L^{2}$ error, and the $L^{2}$ error can also arrive 2.37E-03, which is competitive with the results obtained by virtual
%element method (example 5.1 in \cite{chen2017}). To get the error with same order, we cost about 6.8 seconds, and they cost 191.8 seconds, whcih is enough to show the efficiency of our method. 

%{\color{red}  Although the exact solution is unchanged, the jump and the source term are affected by the coefficients, which require the RNNs to learn different equations. Therefore, adjusting $\gamma_{1}$, $\gamma_{2}$, $\gamma_{3}$ could be a good way to reduce the error when the coefficients differ significantly.}

\begin{table}[!htbp]
    \centering
    \begin{tabular}{|c|c|c|c|c|c|}
    \hline
    $(\beta_{1},\beta_{2})$ & relative $L^{2}$ error & CPU time (s) \\\hline
    ($1,10$)&1.11E-04&6.8 \\\hline
    ($1,10^{3}$)&1.17E-05&6.9 \\\hline
    ($1,10^{5}$)&1.79E-05&6.8  \\\hline
    ($1,10^{7}$)&5.75E-05&6.8  \\\hline
    ($1,10^{-2}$)&7.19E-03&6.7  \\\hline
    ($10^{-2},10^{2}$)&1.08E-05&6.6  \\\hline
    ($10^{-4},10^{4}$)&5.09E-04&6.7 \\\hline
    \end{tabular}
    \caption{Relative $L^2(\Omega)$ errors of $u_\rho$ for different $\beta_{1}$ and $\beta_{2}$, where $m=640,N=10000,r_{1}=2.54$, $r_{2}=0.33$ in Example \ref{ex2}.}\label{tab4}
\end{table}

\begin{example}\label{ex3}
  Two-dimensional Multiple Interfaces Problem. Let $\Omega = (-1,1)^2$, we consider the elliptic interface problem \eqref{model} with three interfaces given as follows:
$$
\Gamma^{0}:x^2+y^2=0.2^{2},
$$
$$
\Gamma^{1}:x^2+y^2=r^{2}(\theta),\quad r(\theta)=0.5-0.1\cos(5\theta),\quad \theta \in [0,2\pi),
$$
$$
\Gamma^{2}:x^2+y^2=0.8^{2}.
$$
The exact solution
$$ 
u(\mathbf{x})=\left\{
\begin{aligned}
&\cos(y)+1.8,\quad \mathbf{x}\in \Omega_{0}, \\
&e^{x}+1.3,\quad \mathbf{x}\in \Omega_{1}, \\
&\sin(x)+0.5,\quad \mathbf{x}\in \Omega_{2}, \\
&-x+\ln(y+2),\quad \mathbf{x}\in \Omega_{3},
\end{aligned}
\right.
$$
$f$, $g^{i}_{1}$, $g^{i}_{2}$, $i=0,1,2$ can be derived from the exact solution and coefficients.
\end{example}

\begin{figure}[!htbp]
      \centering
      \subfloat[approximation $u_\rho$]
      {
        \label{four_1}\includegraphics[width=0.33\textwidth]{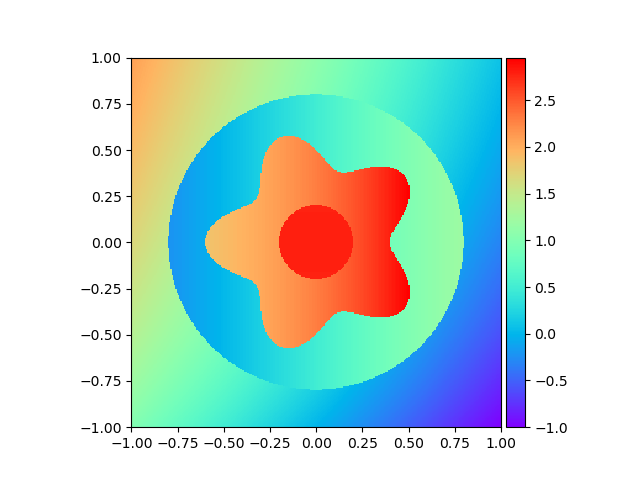}
      }
      \subfloat[real solution $u$]
      {
        \label{four_2}\includegraphics[width=0.33\textwidth]{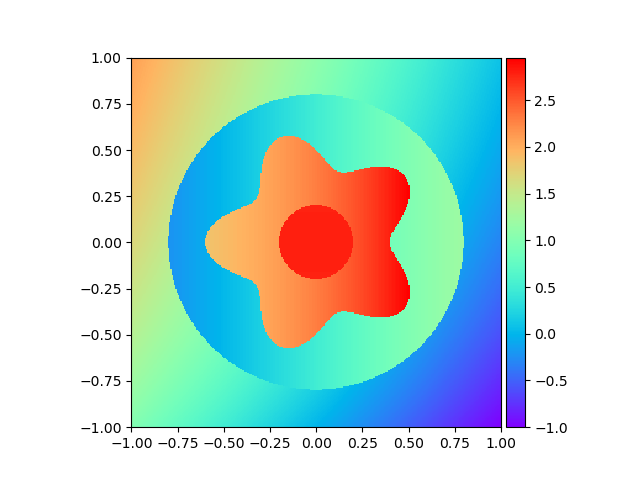}
      }
      \subfloat[absolute error $|u-u_{\rho}|$]
      {
        \label{four_3}\includegraphics[width=0.33\textwidth]{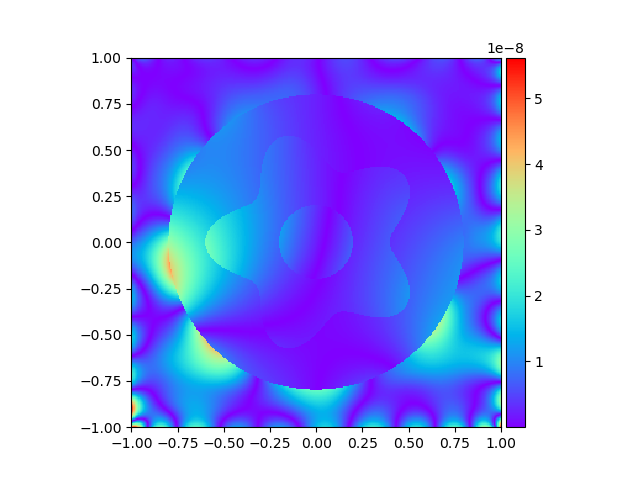}
      }
      \caption{Displays of $u_\rho$, $u$ and $|u-u_{\rho}|$ for $\beta_{i}=i$ ($i=1,2,3,4$) with $m=320$, $N=5000$ in Example \ref{ex3}.}
      \label{figure4}
\end{figure}

We approximate the solution using four RNNs, one for each subdomain. The RNNs have the same structure:  
$n_0 = 2$, $n_1 = m$ and $n_2 = 1$. The weights and biases are randomly drawn from different uniform distributions: $\mathcal{U}(-1.1, 1.1)$, $\mathcal{U}(-0.7, 0.7)$, $\mathcal{U}(-0.3, 0.3)$ and $\mathcal{U}(-1, 1)$. We sample $N$ points from $\Omega$, $\Gamma^{1}$, $\Gamma^{2}$, $\Gamma^{3}$ and $\partial \Omega$ with a ratio of $6:1:1:1:1$.

Figure \ref{figure4} shows the LRNNs solution $u_{\rho}$, the exact solution $u$, and the error $|u-u_{\rho}|$ for $\beta_{1}=1$, $\beta_{2}=2$, $\beta_{3}=3$, $\beta_{4}=4$, $m=320$ and $N=5000$.
Table \ref{tab5} and Table \ref{tab6} report the relative $L^2$ errors for different values of $N$, $m$ and the coefficients $\beta_i$ ($i=1,\cdots,4$). We can see that the LRNNs method performs well for problems with multiple interfaces.
Moreover, our method has a clear advantage over PINNs (Example 4.6 in \cite{hu2020}) in terms of computational efficiency and accuracy under the same settings.

% {\color{green} It should be noted in Table \ref{tab5} is that the error may increase when $N$ is small and $dof(4m)$ is a litter bigger than $N$, the information is not enough to get a more accurate solution}. 

\begin{table}[!htbp]	
	\centering  
	\begin{tabular}{|c|c|c|c|c|c|c|}  
		\hline  
		\diagbox [width=6em] {N}{m}&40&80&160&320&640&1280\cr\cline{1-7}
        500&6.80E-04&3.30E-04&1.68E-03&4.83E-08&2.72E-08&1.59E-08\cr\hline
		    1000&2.50E-04&5.73E-06&6.47E-07&5.58E-06&2.89E-09&1.70E-09\cr\hline
		    2000&1.99E-04&3.74E-06&3.78E-08&9.77E-09&2.19E-07&1.03E-09\cr\hline 
		    5000&2.15E-04&4.74E-06&2.95E-08&4.27E-09&1.30E-09&7.96E-10\cr\hline
        10000&7.93E-04&2.30E-05&7.15E-08&7.57E-09&2.58E-09&7.33E-10\cr\hline
	\end{tabular}  
	\caption{Relative $L^2(\Omega)$ errors with different $N$ and $m$, where $\beta_{i}=i$ ($i=1,2,3,4$) in Example \ref{ex3}.}
	\label{tab5}
\end{table}

\begin{table}[!htbp]
    \centering
    \begin{tabular}{|c|c|c|c|c|c|}
    \hline
    $(\beta_{1},\beta_{2},\beta_{3},\beta_{4})$ & relative $L^{2}$ error & CPU time (s) \\\hline
    ($1,2,3,4$)&3.33E-09&3.5 \\\hline
    ($1,10,10^{2},10^{3}$)&1.89E-07&2.9  \\\hline
    ($1,10^{2},10^{4},10^{6}$)&4.64E-05&3.6  \\\hline
    ($1,10^{-2},10^{-4},10^{-6}$)&7.68E-07&3.6  \\\hline
    ($10^{6},10^{4},10^{2},1$)&7.16E-05&2.9  \\\hline
    ($10^{-6},10^{-4},10^{-2},1$)&4.44E-06&3.5 \\\hline
    \end{tabular}
    \caption{Relative $L^2(\Omega)$ errors of $u_\rho$ for different coefficients, where $m=320$ and $N=5000$ in Example \ref{ex3}.}\label{tab6}
\end{table}

\begin{example}\label{ex6}
High-dimensional Elliptic Interface Problem. Let $\Omega = (0,1)^d$,
we consider an elliptic interface problem \eqref{model} with a hyperplane interface
$$
\Gamma:x_{1}=0.5.
$$
The exact solution is given by
$$ 
u(\mathbf{x})=\left\{
\begin{aligned}
\Vert x \Vert^{2}/d,\quad\mathbf{x}\in \Omega_{1}, \\
\sum\limits_{i=1}^{d}x_{i}/d,\quad\mathbf{x}\in \Omega_{2},
\end{aligned}
\right.
$$
and we can derive $f$, $g_{1}$ and $g_{2}$ from $u$ and the coefficients.
\end{example}

We use two RNNs with the same structure to approximate the solution on each subdomain. The RNNs have two layers with $n_0 = d$, $n_1 = m$ and $n_2 = 1$ respectively. The weights and biases are randomly chosen from the uniform distribution $\mathcal{U}(-0.01,0.01)$. We sample $10^{3},10^{4},2d\times 10^{2}$ points from $\Omega$, $\Gamma$ and $\partial \Omega$, respectively. We use the Monte Carlo method to compute the numerical integration, with $10^{4}$ sample points.

Table \ref{tab13} presents the relative $L^2$ errors for various choices of $d$ and $m$. As in the previous examples, increasing $m$ enhances the accuracy of numerical solutions. The LRNNs method can handle high-dimensional problems effectively, since the degrees of freedom do not need to grow exponentially with the dimension d.

\begin{table}[!htbp]	
	\centering  
	\begin{tabular}{|c|c|c|c|c|}  
		\hline  
		\diagbox [width=6em] {d}{m}&225&450&900&1800\cr\cline{1-5}
        5&8.50E-07&4.74E-07&2.81E-07&1.56E-07\cr\hline
		    10&8.13E-03&3.68E-06&1.79E-06&9.33E-07\cr\hline
		    20&2.92E-02&2.09E-02&1.05E-02&6.98E-05\cr\hline 
	\end{tabular}  
	\caption{Relative $L^2(\Omega)$ errors with different $d$ and $m$, where $\beta_{1}=\beta_{2}=1,r_{1}=r_{2}=0.01$ in Example \ref{ex6}.}
	\label{tab13}
\end{table}

Table \ref{tab14} displays the relative $L^2$  errors and computational time for various choices of diffusion coefficients $\beta_{1}$ and $\beta_{2}$. The results show that the high-dimensional interface problem can be solved accurately and efficiently by the LRNNs method.

\begin{table}[!htbp]
    \centering
    \begin{tabular}{|c|c|c|c|c|c|}
    \hline
    $(\beta_{1},\beta_{2})$ & relative $L^{2}$ error & CPU time (s) \\\hline
    ($1,10^{3}$)&7.24E-05&25.6 \\\hline
    ($1,10^{6}$)&5.38E-04&25.6  \\\hline
    ($1,10^{-3}$)&6.88E-05&26.1  \\\hline
    ($1,10^{-6}$)&5.97E-05&25.3  \\\hline
    ($10^{2},10^{-2}$)&4.54E-03&25.3  \\\hline
    ($10^{-5},10^{5}$)&1.31E-04&25.7 \\\hline
    \end{tabular}
    \caption{Relative $L^2(\Omega)$ errors of $u_\rho$ and CPU time for different $\beta_{1}$ and $\beta_{2}$, where $d=20,m=1800,r_{1}=r_{2}=0.01$ in Example \ref{ex6}.}\label{tab14}
\end{table}

\begin{example}\label{ex4}
   Parabolic Interface Problem with a Fixed Interface. Set $\Omega = (-1,1)^2$ and time interval as $(0,1)$.
We consider the parabolic interface problem \eqref{model_heat} with a circular interface $\Gamma:x^2+y^2=0.5^{2}$. Let 
$$ u(\mathbf{x})=\left\{
\begin{aligned}
&-e^{-t}(8x^{2}+8y^{2}-3.5);(x,y,t)\in \Omega_{1}\times (0,T], \\
&e^{x-t}cos(\frac{\pi}{2}y);(x,y,t)\in \Omega_{2}\times (0,T].
\end{aligned}
\right.
$$
\end{example}

We treat time and space variables equally in the time-dependent problem, so the parabolic interface problem \eqref{model_heat} can be seen as a 3-D interface problem with a cylindrical interface. We approximate the solution using two RNNs, one for each subdomain, with $n_0 = 3$, $n_1 = m$ and $n_2 = 1$. The weights and biases are randomly chosen from $\mathcal{U}(-0.6,0.6)$. We sample $N$ points from $\Omega \times (0,T)$, $\Gamma \times (0,T)$ and $(\partial \Omega \times (0,T))\cup (\Omega \times \{0\})$ with a proportion of $14:3:3$.

\begin{figure}[!htbp]
      \centering     
      \subfloat[approximation $u_\rho$]
      {
        \label{heat01}\includegraphics[width=0.33\textwidth]{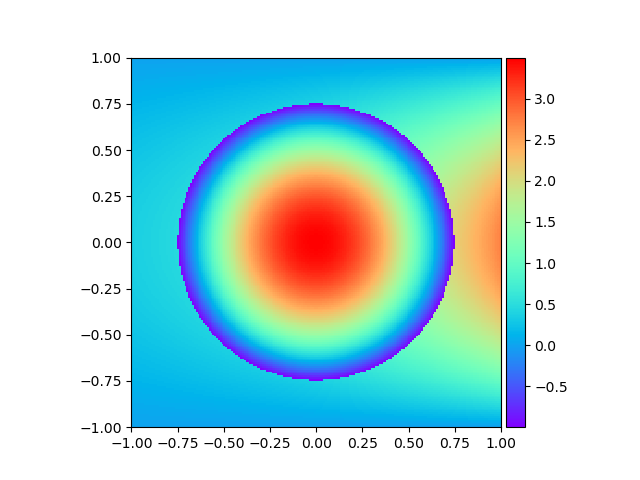}
      }
      \subfloat[real solution $u$]
      {
        \label{heat02}\includegraphics[width=0.33\textwidth]{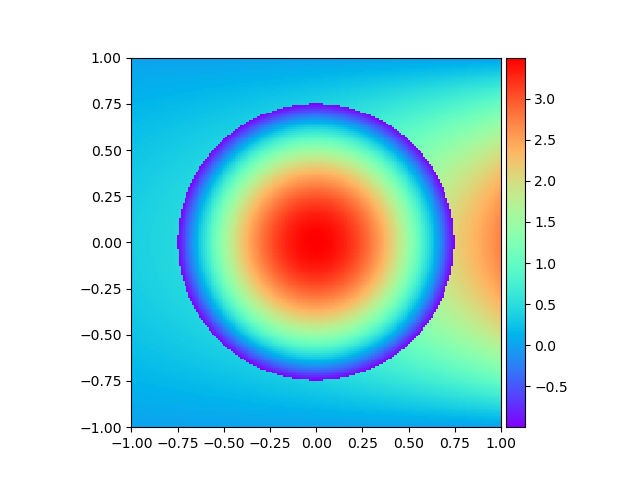}
      }
      \subfloat[absolute error $|u-u_{\rho}|$]
      {
        \label{heat03}\includegraphics[width=0.33\textwidth]{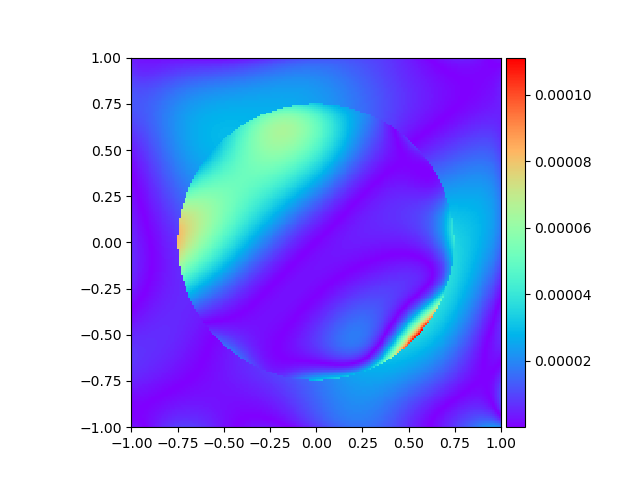}
      }
      \caption{Displays of $u_\rho$, $u$ and $|u-u_{\rho}|$ at $t=0$ in Example \ref{ex4}.}
      \label{figure5}
\end{figure}

\begin{figure}[!htbp]
    \centering     
    \subfloat[approximation $u_\rho$]
    {
      \label{heat11}\includegraphics[width=0.33\textwidth]{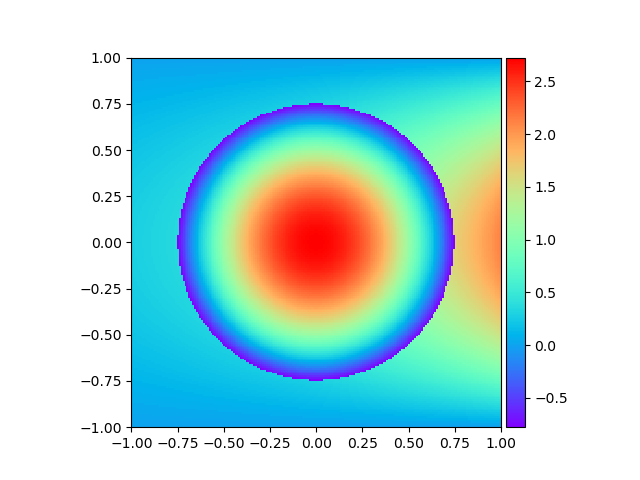}
    }
    \subfloat[real solution $u$]
    {
      \label{heat12}\includegraphics[width=0.33\textwidth]{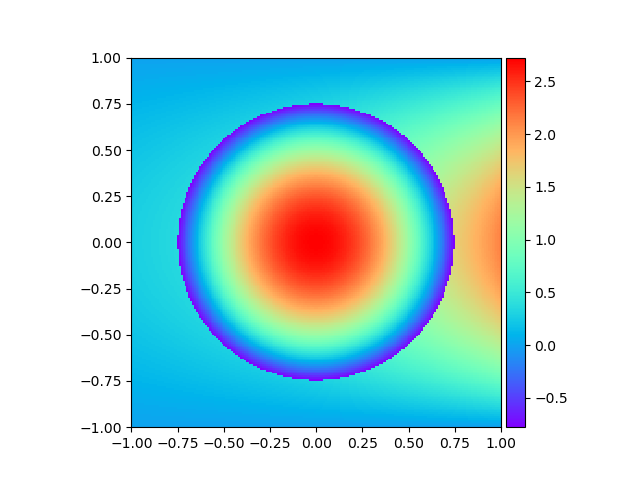}
    }
    \subfloat[absolute error $|u-u_{\rho}|$]
    {
      \label{heat13}\includegraphics[width=0.33\textwidth]{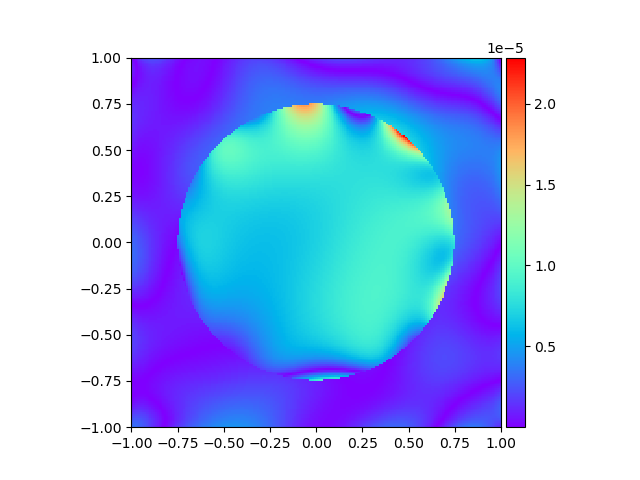}
    }
    \caption{Displays of $u_\rho$, $u$ and $|u-u_{\rho}|$ at $t=0.25$ in Example \ref{ex4}.}
    \label{figure6}
\end{figure}

\begin{figure}[!htbp]
    \centering     
    \subfloat[approximation $u_\rho$]
    {
      \label{heat21}\includegraphics[width=0.33\textwidth]{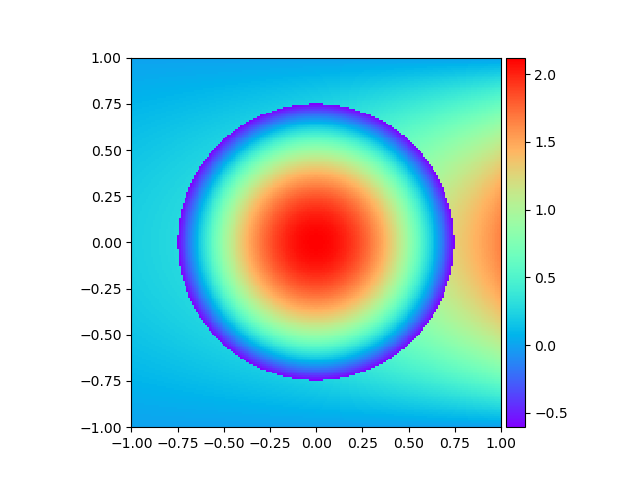}
    }
    \subfloat[real solution $u$]
    {
      \label{heat22}\includegraphics[width=0.33\textwidth]{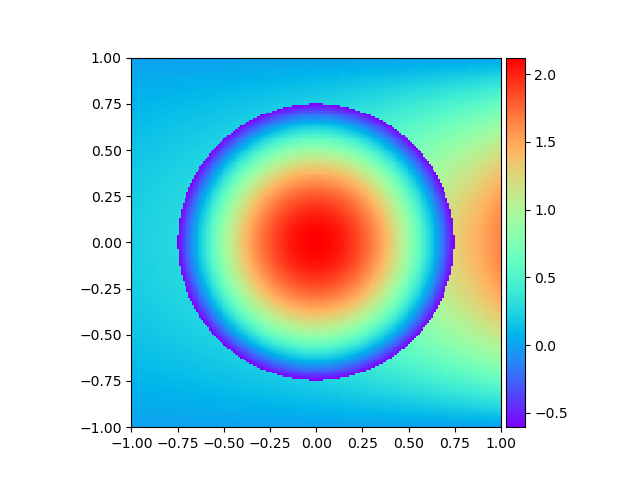}
    }
    \subfloat[absolute error $|u-u_{\rho}|$]
    {
      \label{heat23}\includegraphics[width=0.33\textwidth]{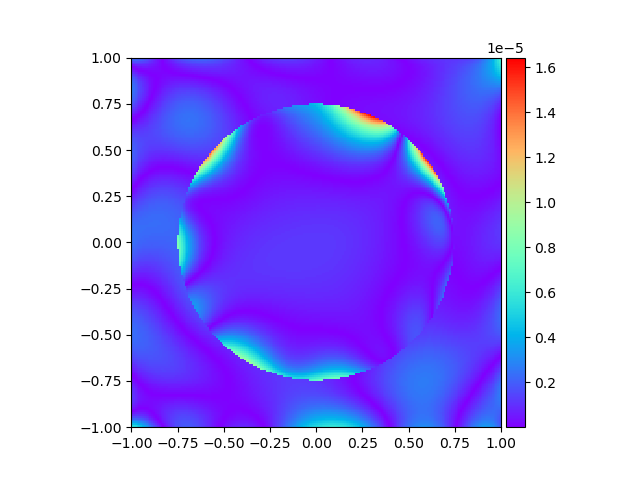}
    }
    \caption{Displays of $u_\rho$, $u$ and $|u-u_{\rho}|$ at $t=0.5$ in Example \ref{ex4}.}
    \label{figure7}
\end{figure}

\begin{figure}[!htbp]
  \centering     
  \subfloat[approximation $u_\rho$]
  {
    \label{heat31}\includegraphics[width=0.33\textwidth]{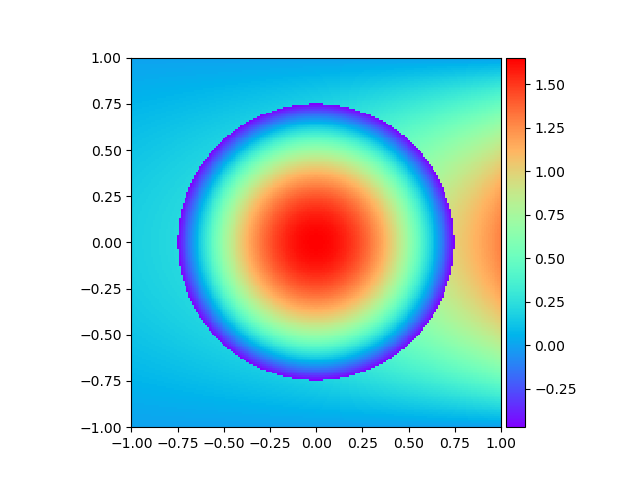}
  }
  \subfloat[real solution $u$]
  {
    \label{heat32}\includegraphics[width=0.33\textwidth]{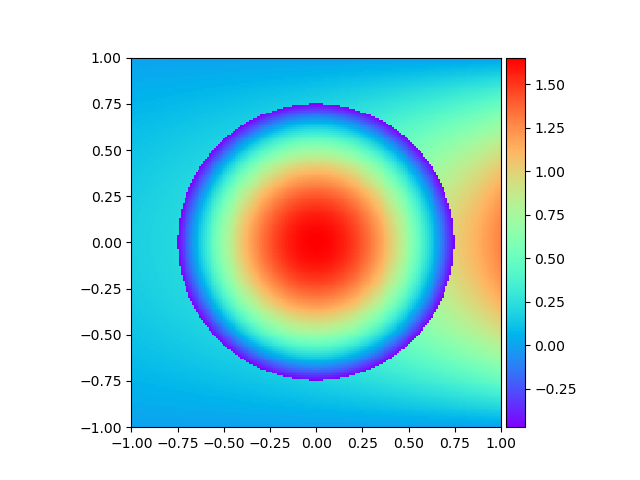}
  }
  \subfloat[absolute error $|u-u_{\rho}|$]
  {
    \label{heat33}\includegraphics[width=0.33\textwidth]{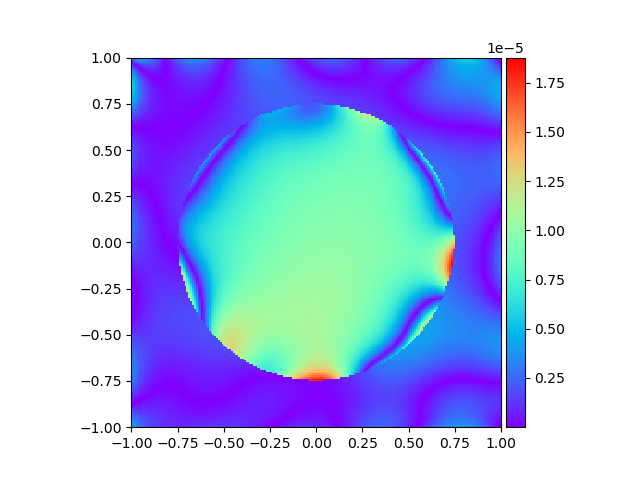}
  }
  \caption{Displays of $u_\rho$, $u$ and $|u-u_{\rho}|$ at $t=0.75$ in Example \ref{ex4}.}
  \label{figure8}
\end{figure}

\begin{figure}[!htbp]
  \centering     
  \subfloat[approximation $u_\rho$]
  {
    \label{heat41}\includegraphics[width=0.33\textwidth]{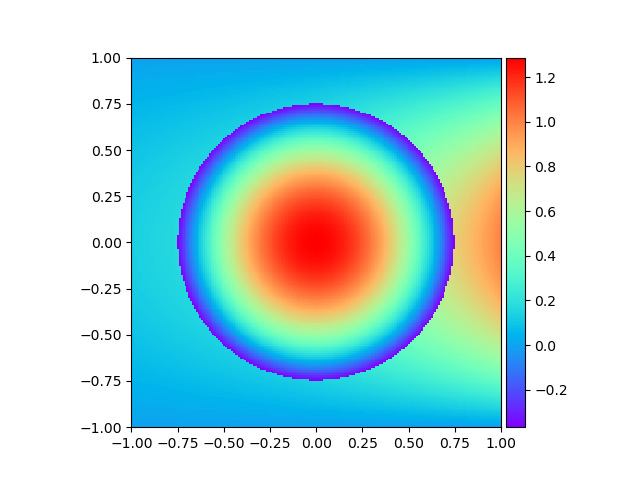}
  }
  \subfloat[real solution $u$]
  {
    \label{heat42}\includegraphics[width=0.33\textwidth]{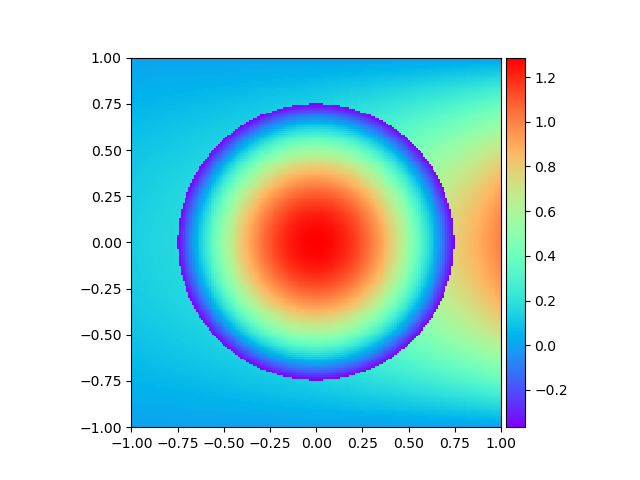}
  }
  \subfloat[absolute error $|u-u_{\rho}|$]
  {
    \label{heat43}\includegraphics[width=0.33\textwidth]{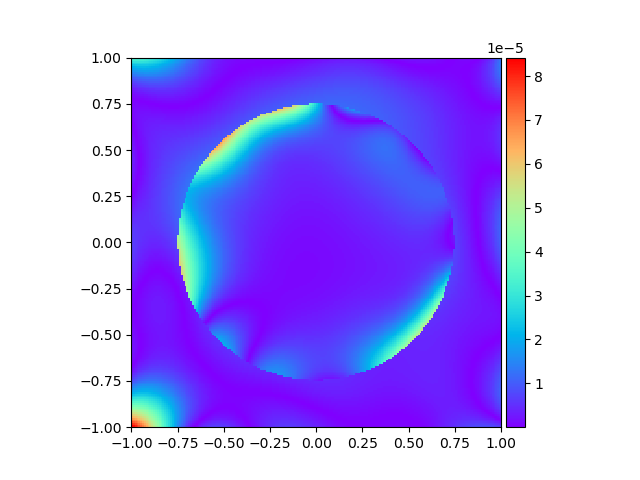}
  }
  \caption{Displays of $u_\rho$, $u$ and $|u-u_{\rho}|$ at $t=1$ in Example \ref{ex4}.}
  \label{figure9}
\end{figure}

Figure \ref{figure5}--\ref{figure9} show the RNN solution $u_\rho$, the exact solution $u$, and the error $|u-u_{\rho}|$ at different time: $t=0, 0.25, 0.5, 0.75, 1$. Here, we choose $m=320$, $N=5000$ for case of $\beta_{1}=1$, $\beta_{2}=1$. We observe that the LRNNs method achieves high accuracy with no error accumulation over time, unlike traditional numerical methods. 

Table \ref{tab7} reports the relative $L^2$ errors for different values of $N$ and $m$. Table \ref{tab8} shows the relative $L^2$ errors and the computational time for different values of the coefficients. This time-dependent problem does not require any time discretization by the space-time LRNNs method, and only involves solving a linear least square problem. These results demonstrate the robustness and effectiveness of the LRNNs method for the parabolic interface problem.

\begin{table}[!htbp]	
	\centering  
	\begin{tabular}{|c|c|c|c|c|c|}  
		\hline  
		\diagbox [width=6em] {N}{m}&80&160&320&640&1280\cr\cline{1-6}
		    1000&4.05E-03&3.70E-04&1.01E-04&1.57E-05&1.44E-05\cr\hline 
		    2000&3.37E-03&2.17E-04&1.86E-05&1.09E-06&7.04E-07\cr\hline
        5000&3.37E-03&1.52E-04&5.94E-06&4.49E-07&7.39E-08\cr\hline
        10000&2.48E-03&1.87E-04&4.99E-06&3.34E-07&7.84E-08\cr\hline 
        20000&3.29E-03&1.58E-04&5.37E-06&3.27E-07&1.00E-07\cr\hline 
	\end{tabular}  
	\caption{Relative $L^2(\Omega)$ errors with different $N$ and $m$, where $\beta_{1}=\beta_{2}=1$, $r_{1}=r_{2}=0.6$ in Example \ref{ex4}. }
	\label{tab7}
\end{table}

\begin{table}[!htbp]
  \centering
  \begin{tabular}{|c|c|c|c|c|c|}
  \hline
  $(\beta_{1},\beta_{2})$ & relative $L^{2}$ error & CPU time (s) \\\hline
  ($1,10^{2}$)&1.68E-06&5.1 \\\hline
  ($1,10^{5}$)&4.14E-04&5.3  \\\hline
  ($10^{-2},10^{2}$)&3.50E-04&5.1  \\\hline
  ($10^{2},10^{-2}$)&3.95E-05&5.7  \\\hline
  ($10^{4},1$)&1.40E-04&5.5  \\\hline
  ($10^{2},1$)&7.89E-07&5.4 \\\hline
  \end{tabular}
  \caption{Relative $L^2(\Omega)$ errors of $u_\rho$ for different $\beta_{1}$ and $\beta_{2}$, where $m=320$ and $N=5000$ in Example \ref{ex4}.}\label{tab8}
\end{table}

\begin{example}\label{ex5}
Parabolic Interface Problem with a Dynamic Interface. Set $\Omega = (-1,1)^2$ and time interval as $(0,1)$.
We consider the parabolic interface problem \ref{model_heat} with a circular interface $\Gamma:x^2+y^2=(0.5+0.3t)^{2}$. The exact solution is
$$ u(\mathbf{x})=\left\{
\begin{aligned}
&-e^{-t}(8x^{2}+8y^{2}-3.5);(x,y,t)\in \Omega_{1}\times (0,T], \\
&e^{x-t}cos(\frac{\pi}{2}y);(x,y,t)\in \Omega_{2}\times (0,T].
\end{aligned}
\right.
$$
\end{example}

In this example, we consider a dynamic interface problem, where the interface is shaped like: the curved surface of a cone’s frustum in the space-time domain. We use two RNNs with $n_0 = 3$, $n_1 = m$ and $n_2 = 1$ to approximate the solution on each subdomain. The weights and biases are randomly chosen from $(-1,1)$. We sample $N$ points from $\Omega \times (0,T)$, $\Gamma \times (0,T)$ and $(\partial \Omega \times (0,T))\cup (\Omega \times \{0\})$ with a ratio of $14:3:3$. 

\begin{figure}[!htbp]
     \centering     
     \subfloat[approximation $u_\rho$]
     {
       \label{heat51}\includegraphics[width=0.33\textwidth]{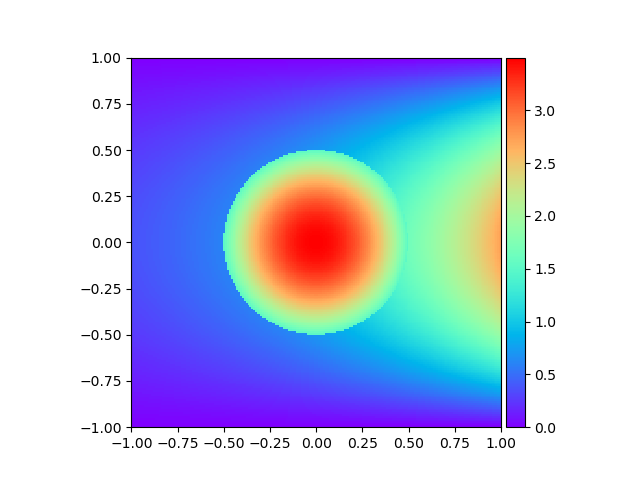}
     }
     \subfloat[real solution $u$]
     {
       \label{heat52}\includegraphics[width=0.33\textwidth]{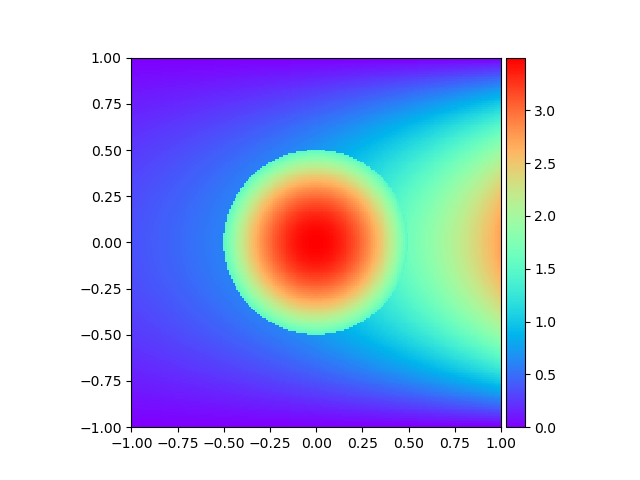}
     }
     \subfloat[absolute error $|u-u_{\rho}|$]
     {
       \label{heat53}\includegraphics[width=0.33\textwidth]{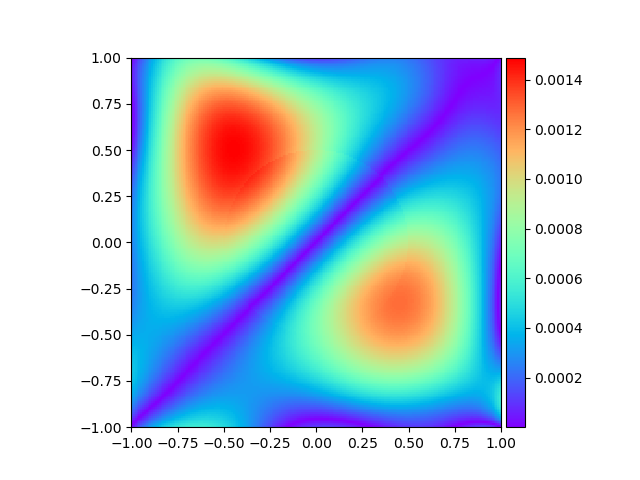}
     }
     \caption{Displays of $u_\rho$, $u$ and $|u-u_{\rho}|$ at $t=0$ in Example \ref{ex5}.}
     \label{figure11}
\end{figure}

\begin{figure}[!htbp]
   \centering     
   \subfloat[approximation $u_\rho$]
   {
     \label{heat61}\includegraphics[width=0.33\textwidth]{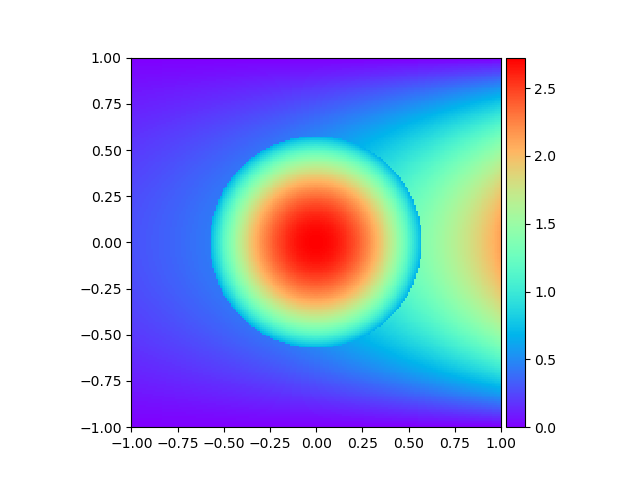}
   }
   \subfloat[real solution $u$]
   {
     \label{heat62}\includegraphics[width=0.33\textwidth]{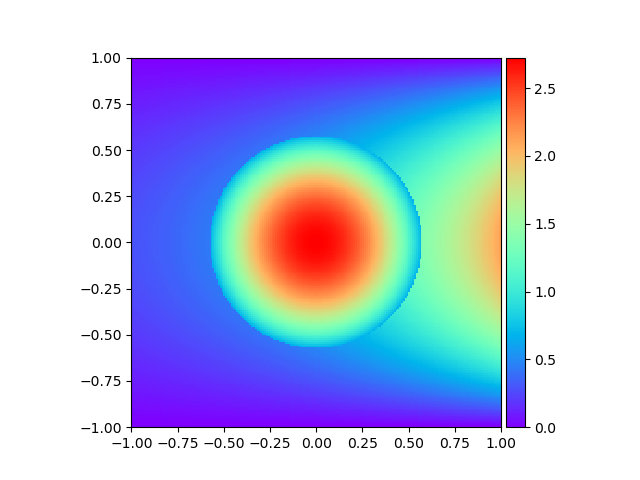}
   }
   \subfloat[absolute error $|u-u_{\rho}|$]
   {
     \label{heat63}\includegraphics[width=0.33\textwidth]{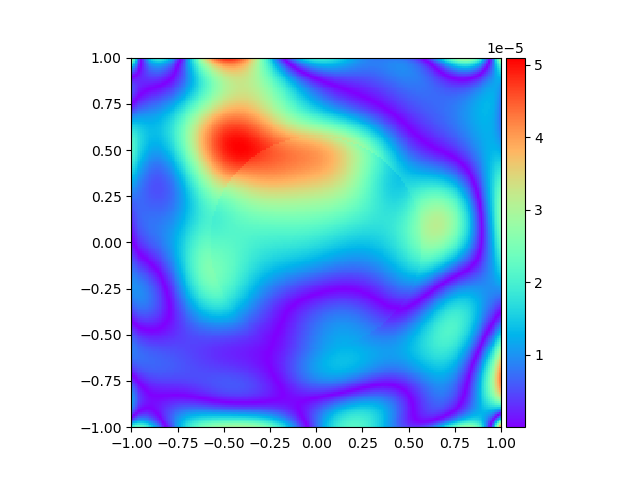}
   }
   \caption{Displays of $u_\rho$, $u$ and $|u-u_{\rho}|$ at $t=0.25$ in Example \ref{ex5}.}
   \label{figure12}
\end{figure}

\begin{figure}[!htbp]
   \centering     
   \subfloat[approximation $u_\rho$]
   {
     \label{heat71}\includegraphics[width=0.33\textwidth]{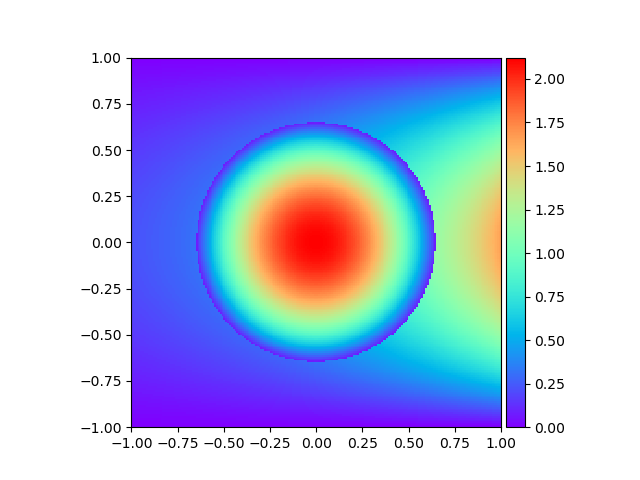}
   }
   \subfloat[real solution $u$]
   {
     \label{heat72}\includegraphics[width=0.33\textwidth]{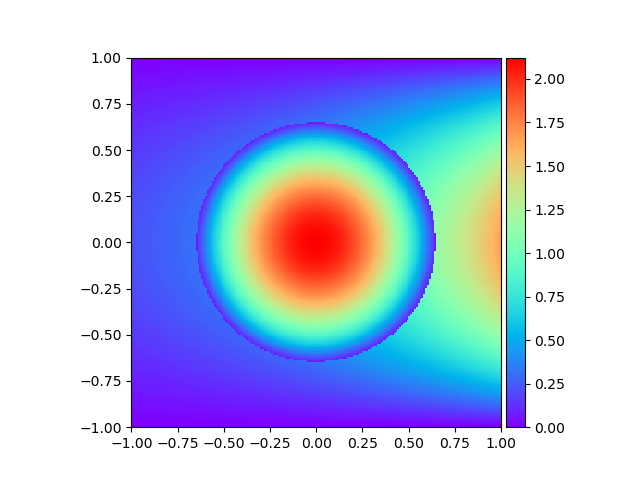}
   }
   \subfloat[absolute error $|u-u_{\rho}|$]
   {
     \label{heat73}\includegraphics[width=0.33\textwidth]{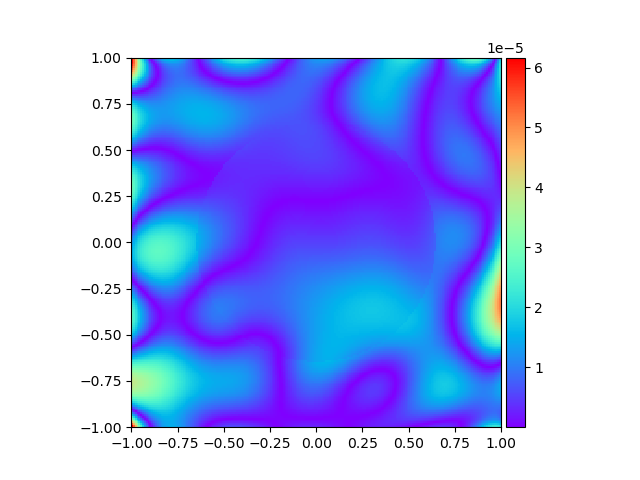}
   }
   \caption{Displays of $u_\rho$, $u$ and $|u-u_{\rho}|$ at $t=0.5$ in Example \ref{ex5}.}
   \label{figure13}
\end{figure}

\begin{figure}[!htbp]
 \centering     
 \subfloat[approximation $u_\rho$]
 {
   \label{heat81}\includegraphics[width=0.33\textwidth]{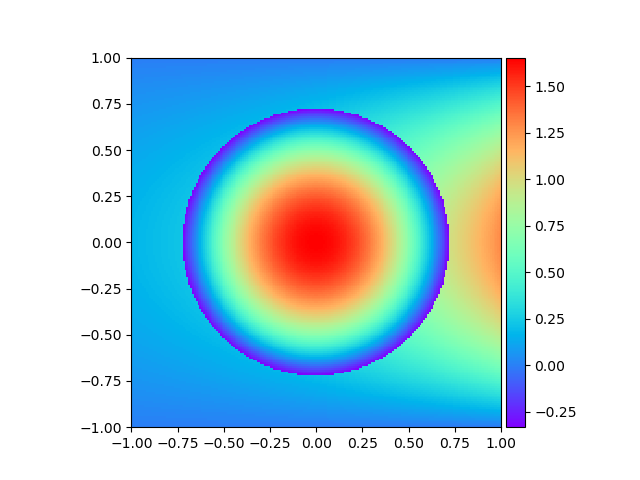}
 }
 \subfloat[real solution $u$]
 {
   \label{heat82}\includegraphics[width=0.33\textwidth]{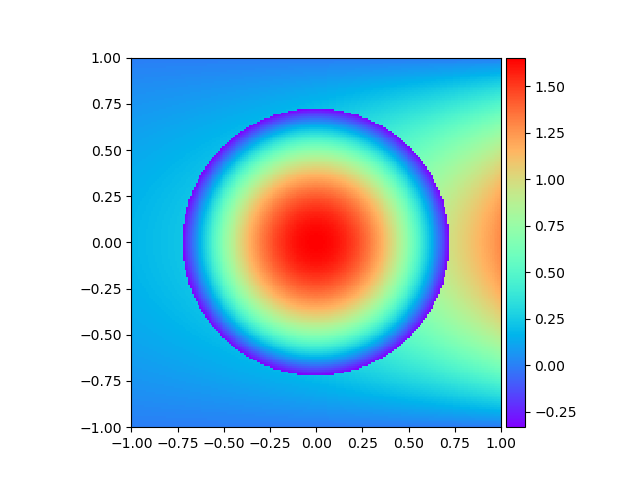}
 }
 \subfloat[absolute error $|u-u_{\rho}|$]
 {
   \label{heat83}\includegraphics[width=0.33\textwidth]{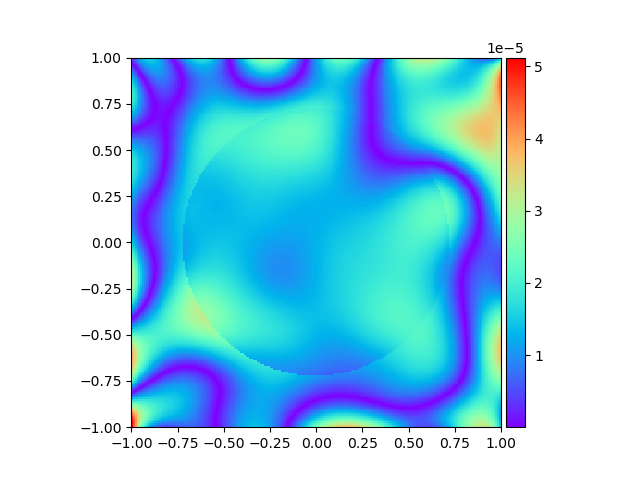}
 }
 \caption{Displays of $u_\rho$, $u$ and $|u-u_{\rho}|$ at $t=0.75$ in Example \ref{ex5}.}
 \label{figure14}
\end{figure}

\begin{figure}[!htbp]
 \centering     
 \subfloat[approximation $u_\rho$]
 {
   \label{heat91}\includegraphics[width=0.33\textwidth]{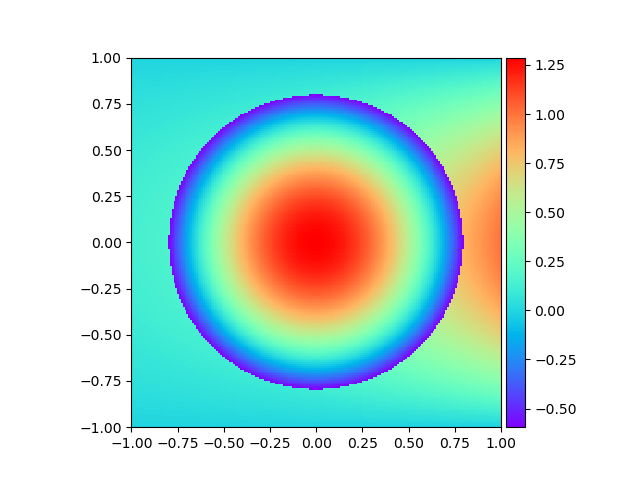}
 }
 \subfloat[real solution $u$]
 {
   \label{heat92}\includegraphics[width=0.33\textwidth]{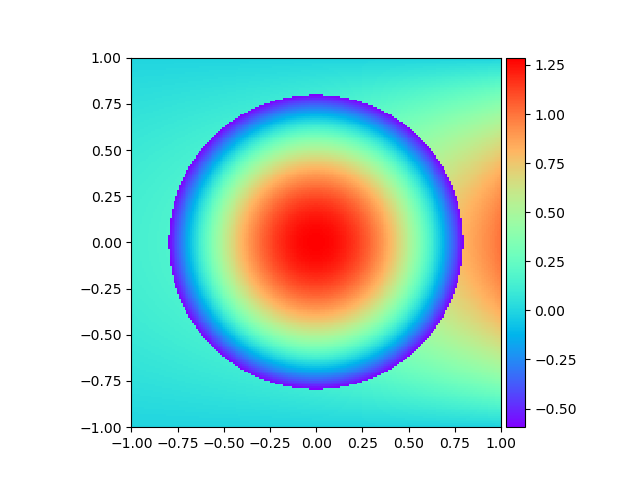}
 }
 \subfloat[absolute error $|u-u_{\rho}|$]
 {
   \label{heat93}\includegraphics[width=0.33\textwidth]{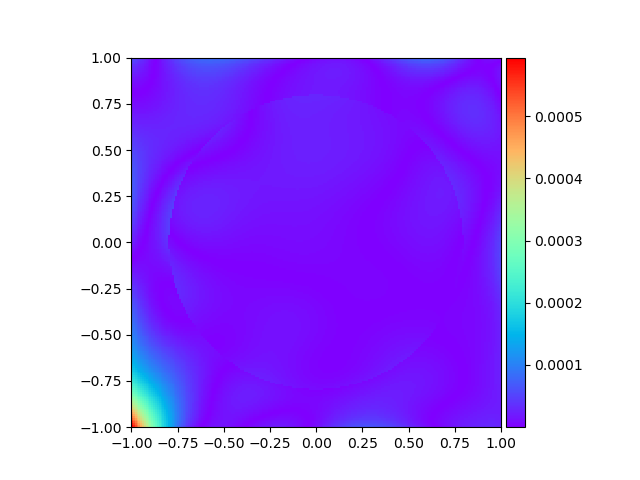}
 }
 \caption{Displays of $u_\rho$, $u$ and $|u-u_{\rho}|$ at $t=1$ in Example \ref{ex5}.}
 \label{figure15}
\end{figure}

Figures \ref{figure11} to \ref{figure15} illustrate the comparison between the exact solution u, the RNN solution $u_\rho$, and the absolute error $|u-u_{\rho}|$ at various time points: $t=0, 0.25, 0.5, 0.75, 1$. We use $\beta_{1}=1$, $\beta_{2}=1$, $m=320$ and $N=5000$ for this case. The interface evolves over time, and the LRNNs method preserves high accuracy without accumulating errors.

In Table \ref{tab11}, we present the relative $L^2$ errors for different choices of $N$ and $m$. We also investigate the influence of the coefficients on the relative $L^2$ errors in Table \ref{tab12}. These results confirm the robustness and effectiveness of the space-time LRNNs method for solving the parabolic moving interface problem.

\begin{table}[!htbp]	
 \centering  
 \begin{tabular}{|c|c|c|c|c|c|}  
   \hline  
   \diagbox [width=6em] {N}{m}&80&160&320&640&1280\cr\cline{1-6}
       1000&8.10E-03&2.33E-03&2.37E-04&7.18E-05&5.78E-05\cr\hline 
       2000&6.80E-03&6.97E-04&8.23E-05&7.37E-06&1.63E-06\cr\hline
       5000&6.77E-03&5.72E-04&2.97E-05&2.22E-06&2.48E-07\cr\hline
       10000&6.12E-03&5.48E-04&2.24E-05&1.28E-06&2.27E-07\cr\hline 
       20000&6.08E-03&4.96E-04&2.77E-05&2.56E-06&2.07E-07\cr\hline 
 \end{tabular}  
 \caption{Relative $L^2(\Omega)$ errors with different $N$ and $m$, where $\beta_{1}=\beta_{2}=1$, $r_{1}=1$, $r_{2}=1$ in Example \ref{ex5}. }
 \label{tab11}
\end{table}

\begin{table}[!htbp]
 \centering
 \begin{tabular}{|c|c|c|c|c|c|}
 \hline
 $(\beta_{1},\beta_{2})$ & relative $L^{2}$ error & CPU time (s) \\\hline
 ($1,10^{2}$)&8.75E-06&3.7 \\\hline
 ($1,10^{5}$)&4.34E-03&3.7  \\\hline
 ($10^{-2},10^{2}$)&1.30E-03&4.4  \\\hline
 ($10^{2},10^{-2}$)&5.20E-06&3.7  \\\hline
 ($10^{4},1$)&8.35E-05&3.7  \\\hline
 ($10^{2},1$)&9.02E-07&4.3 \\\hline
 \end{tabular}
 \caption{Relative $L^2(\Omega)$ errors of $u_\rho$ for different $\beta_{1}$ and $\beta_{2}$, where $m=320$ and $N=5000$ in Example \ref{ex5}.}\label{tab12}
\end{table}

\section{Summary}
\label{sec5}

We develop LRNNs methods to solve interface problems with multiple RNNs. Each RNN works on a different sub-domain and approximates the solution there. We use the least-square method to solve the linear system. The LRNNs method has several benefits: it does not need a mesh, so it can handle complex interfaces easily; it can achieve high accuracy with fewer degrees of freedom; it does not need to solve an optimization problem, which saves a lot of computation time; it can solve time-dependent problems with a space-time approach, which avoids error accumulation from time iteration; it is robust for different coefficients.

We plan to explore more features and applications of the LRNNs method in the future. Some of our research directions are: finding a suitable range of weights (and bias) for the RNNs, which may depend on the source term and other information; developing an adaptive strategy for sampling points to improve the efficiency of the LRNNs method; establishing the numerical analysis for the LRNNs methods; applying the LRNNs method to other complex problems, such as two-phase flow problems and fluid-structure interaction problems.

\end{document}